\documentclass[a4paper,12pt]{article}

\ifx\TeXtures\undefined
\usepackage[psamsfonts]{amsfonts}%
\usepackage{amsmath}%
\else
\usepackage{amsmath}%
\usepackage{amsfonts}%
\usepackage{times}
\fi
\usepackage{amsthm} 
\usepackage{epsfig}         
\newcommand{\DATUM}{30-11-2007}              
\newcommand{\change}
{{\marginpar{\#}}}        
\newcommand{\comma}{\: ,}              
\newcommand{\period}{\: .}             
\newcommand{\cA}{{\cal A}}

\newcommand{\cE}{{\cal E}}
\newcommand{\cF}{{\cal F}}

\newcommand{\cH}{{\cal H}}

\newcommand{\cJ}{{\cal J}}
\newcommand{\cK}{{\cal K}}
\newcommand{\cL}{{\cal L}}         

\newcommand{\cS}{{\cal S}}
\newcommand{\cT}{{\cal T}}

\newcommand{\field}[1]{\mathbb{#1}}
\newcommand{\R}{\field{R}}            
\newcommand{\Z}{\field{Z}}            
\newcommand{\N}{\field{N}}            
\newcommand{\C}{\field{C}}            

\newcommand{\rL}{{\rm L}}                 
\newcommand{\rS}{{\rm S}}                 
\newcommand{\rB}{{\rm B}}                 
\newcommand{\rH}{{\rm H}}

\newcommand{\ph}{{\rm x}^\ast}
\newcommand{\phz}{{\rm z}^\ast}
\newcommand{\pht}{{\rm t}^\ast}

\newcommand{\rRe}{{\rm Re}}               

\newcommand{\support}{{\rm supp}\, }

\newcommand{\cirS}{\mathop{\bigcirc\kern -.73em {\scriptstyle{\rm S}}}}

\newcommand{\qmbox}[1]{\quad\mbox{#1}\quad}
\newtheorem{theorem}{Theorem}[section]                
\newtheorem{lemma}[theorem]{Lemma}             
\newtheorem{corollary}[theorem]{Corollary}     
\newtheorem{remark}[theorem]{Remark}           
\newtheorem{proposition}[theorem]{Proposition} 
\newtheorem{example}[theorem]{Example}
\theoremstyle{plain}
\newcommand{\tend}{\rightarrow}
\newcommand{\donne}{\mapsto}
\newcommand{\dans}{\longrightarrow}

\newcommand{\impl}{\Longrightarrow}

\newcommand{\un}{1\hskip-0.6ex {\sf I}}
\newcommand{\Pf}{\vspace*{-2mm}{\bf Proof:}\, }
\newcommand{\Pfof}[1]{{\bf Proof of #1:}\, }

\newcommand{\rond}{\!\circ\!}

\newcommand{\barre}[1]{\overline{#1}}

\newcommand {\tr}{{\rm tr}}
\renewcommand {\l}{\left}
\newcommand {\ri}{\right}
\newcommand{\J}{{\rm coll}}
\newcommand{\ar}{\rightarrow}
\newcommand {\bH}{{\mathbb H}}
\newcommand{\beq}{\begin{equation}}
\newcommand{\Leq}[1]{\label{#1}\end{equation}}
\newcommand{\eeq}{\end{equation}}
\newcommand{\beqno}{\begin{eqnarray*}}
\newcommand{\eeqno}{\end{eqnarray*}}
\newcommand{\bem}{\l(\! \begin{array}}
\newcommand{\eem}{\end{array}\!\ri)}
\newcommand{\bsm}{\left(\begin{smallmatrix}} 
\newcommand{\esm}{\end{smallmatrix}\right)}  
\newcommand {\eh}{{\textstyle \frac{1}{2}}}
\newcommand {\LA}{\left\langle}
\newcommand {\RA}{\right\rangle}
\newcommand{\idty}{{\rm 1\mskip-4mu l}} 

\begin{document}

\setcounter{section}{0}

\title{Semiclassical resolvent estimates for Schr\"odinger operators with 
Coulomb singularities.}
\author{Fran\c cois Castella
\thanks{IRMAR \& IRISA - Universit\'{e} de Rennes 1, 
Campus de Beaulieu, 35042, Rennes Cedex, France. 
e-mail :  francois.castella@univ-rennes1.fr, web: http://perso.univ-rennes1.fr/francois.castella/} 
\and
Thierry Jecko
\thanks{IRMAR - Universit\'{e} de Rennes 1, 
Campus de Beaulieu, 35042, Rennes Cedex, France. 
e-mail : jecko@univ-rennes1.fr, web: http://perso.univ-rennes1.fr/thierry.jecko/}
\and
Andreas Knauf
\thanks{Mathematisches Institut, 
Universit\"at Erlangen-N\"urnberg,
Bismarckstr.\ 1 1/2, D--91054 Erlangen, Germany. 
e-mail: knauf@mi.uni-erlangen.de, web: www.mi.uni-erlangen.de/$\sim$knauf}
}
\date{\DATUM}
\maketitle

\begin{abstract}
Consider the Schr\"odinger operator with
semiclassical parameter $h$, in the limit where $h$ goes to zero.
When the involved long-range potential is smooth,
it is well known that the boundary values of the operator's resolvent 
at a positive energy $\lambda$ are bounded by $O(h^{-1})$
if and only if the associated Hamilton 
flow is non-trapping at energy $\lambda$. In the present paper,
we extend this result 
to the case where the potential may possess Coulomb singularities. Since the 
Hamilton flow then is not complete in general, our analysis 
requires the use of an appropriate regularization. 

\vspace{2mm}

\noindent
{\bf Keywords:} Schr\"odinger operators, Coulomb singularities, 
regularization, semiclassical resolvent estimates, resonances,
semiclassical measure.
\end{abstract}

\newpage


{\footnotesize \tableofcontents}



\section{Introduction.}
\label{intro}
\setcounter{equation}{0}

In the late eighties and the beginning of the nineties, many semiclassical 
results were obtained in stationary scattering theory. In this setting, the 
long time evolution of a system is studied via the resolvent, which appears 
in representation formulae for the main scattering objects. 
One can distinguish two complementary domains: on the one hand semiclassical 
results concerning scattering objects at non-trapping energies  
(when resonances are negligible), and on the other hand studies of
resonances and of their influence on scattering objects. 
We refer to \cite{gm2,hs,kmw,ma2,rt} and also to \cite{r2} for an overview of the subject. 
These results often show a Bohr correspondence principle for the scattering states. \\
Many studies treat (non-relativistic) molecular systems described by a (many body) 
Schr\"o\-din\-ger operator. From a physical point of view, it is natural to let
the potential admit Coulomb singularities in that context. In the spectral
analysis of the 
operator, these singularities do not produce difficulties in 
dimension $\geq 3$, thanks to Hardy's inequality (cf. (\ref{hardy})). \\
In the semiclassical regime however, little is known when Coulomb
singularities occur. We point out  
the propagation results in \cite{gk,ke,kn}. In the above mentioned domains 
of stationary scattering theory, we do not know of any semiclassical result, 
except that of \cite{kmw,w3}. We think that the main obstacle stems from the difficulty to 
develop a semiclassical version of Mourre's theory (cf. \cite{gm,m,rt}) in this situation.
This task is performed in \cite{kmw} when all singularities are repulsive, a situation 
where the associated classical Hamilton flow is complete. 
Recently semiclassical resolvent estimates (and further interesting results) were obtained by Wang
in \cite{w3} but in a non optimal framework (see comments below). When attractive singularities 
occur, the classical flow is not complete anymore, while it can be regularized
(cf. \cite{gk,ke,kn}).\\
The aim of this article is to contribute to the development of such a
semiclassical analysis of molecular scattering.
In \cite{jec0,jec1}, the author faced similar difficulties in the 
study of a matricial Schr\"odinger operator. 
He adapted in \cite{jec2,jec3} an alternative approach, previously used in \cite{b}. 
We here follow the same approach, combined with ideas from \cite{cj,gk,kn,w}, 
in order to extend, in the case of potentials with arbitrary
Coulomb singularities, a result established
in \cite{kmw,rt}. \\
We now introduce some notation and present the main results of this paper. 
\subsection{The Schr\"odinger operators.}
\label{1.1}
Let $d\in \N:=\{0,1,2,\ldots \}$ with $d\geq 2$. For $x\in \R^d$, we denote by 
$|x|$ the usual norm of $x$ and we set $\langle x\rangle :=(1+|x|^2)^{1/2}$. We 
denote by $\Delta_x$ the Laplacian in $\R^d$. We consider a long-range potential 
$V$ which is smooth except at $N'$ Coulomb singularities ($N'\in \N^\ast$) 
located at the sites $s_j$ where $j\in\{1,2,\ldots ,N'\}$. 
Let 
\[\hat{M}:=\R^d\setminus \cS\qmbox{, with} \cS:=\{s_j; 1\leq j\leq N'\}\]
and $R_0:=\max \{|s_j|+1; 
1\leq j\leq N'\}$. Technically, we take $V\in C^\infty (\hat{M};\R)$ such that 
\begin{equation}\label{decroit}
\exists \rho>0\, ;\, \forall \alpha \in \N^d\comma \, \forall x\in \R^d\comma |x|>R_0
\comma  \hspace{.4cm}|\partial _x^\alpha V(x)|\ = \ O_\alpha \bigl(\langle x\rangle 
^{-\rho -|\alpha|}\bigr)
\period 
\end{equation}
Furthermore, we assume that for all $j\in \{1,2,\ldots ,N'\}$, we can find smooth functions $f_j, W_j$ in $C_0^\infty(\R^d ;\R)$ such that $f_j(s_j)\neq 0$ and, near $s_j$, 
\begin{equation}\label{sing-sj}
V(x)\ = \ \frac{f_j(x)}{|x-s_j|}\, +\, W_j(x)\period 
\end{equation}
If $f_j(s_j)<0$ (resp. $f_j(s_j)>0$), we say that $s_j$ is an {\em
attractive} (resp.\ a {\em repulsive}) 
{\em Coulomb singularity}. Let $N\geq 0$ be the number of attractive singularities. We may assume that 
they are labelled by $\{1,2,\ldots ,N\}$. \\
Given some $h_\ast\in ]0;1[$, we introduce a semiclassical parameter 
$h\in ]0;h_\ast]$.
The semiclassical Schr\"odinger operator is given by 
$P(h):=-h^2\Delta_x + V$, acting 
in $\rL^2(\R^d)$.
\\
Under the previous assumptions, it is well known that $P(h)$ is 
self-adjoint (see \cite{cfks,k,rs2}). When $d\geq 3$, this fact follows from Hardy's inequality
(cf. (\ref{hardy})) and from Kato's theorem on relative boundedness.
The domain of $P(h)$ then is the Sobolev space 
$\rH ^2(\R^d)$, i.e. the domain of the Laplacian.
When $d=2$, self-adjointness  follows when considering the 
quadratic form associated with $P(h)$ and using Kato's theorem on relative boundedness for 
forms: the form is seen to be closable and bounded below, and the associated
self-adjoint operator is $P(h)$ \cite{ch,kn3}.
The situation is rather different for $d=1$ 
(see Section I.1 in \cite{rs2}), which is the reason why we exclude
this dimension  here. 
\subsection{The function spaces and main notation.}

For $z$ belonging to the resolvent set $\rho (P(h))$ of $P(h)$, we set $R(z;h):=(P(h)-z)^{-1}$.
We are interested in the size of the resolvent $R(z;h)$ as a bounded 
operator from some space $\rS$ into its dual $\rS^\ast$, i.e.\ as an 
element of the space $\cL(\rS ;\rS^\ast)$.  We denote by $\|\cdot \|_{\rS ,\rS^\ast}$ 
the usual operator norm on $\cL(\rS ;\rS^\ast)$. If $\rS =\rL^2(\R^d)$, 
we also use the notation $\|\cdot \|$ in place of $\|\cdot \|_{\rS ,\rS^\ast}$. 
The relevant spaces $\rS$ are introduced below. \\
If $a$ is a measurable subset of $\R^d$, we denote by $\|\cdot\|_a$ (resp. 
$\langle \cdot ,\cdot \rangle_a$) the usual norm (resp. the 
right linear scalar product) of $\rL^2(a)$ (and we skip the subscript $a$ if $a=\R^d$). 
For $s\in \R$, we denote by $\rL _s^2$ the weighted 
$\rL^2$-space of measurable functions $f$ such that $x\donne 
\langle x\rangle^sf(x)$ belongs to $\rL^2(\R^d)$. Its dual space is identified 
with $\rL _{-s}^2$. For $j\in \Z$, we set 
\beq
c_j:=\{x\in \R^d; 2^{j-1}< |x|\leq 2^j\}\qmbox{and} 
c=\{x\in \R^d; |x|\leq 1\}.
\Leq{c:c} 
Let $\rB$ (resp. its homogeneous 
version $\dot{\rB}$) be the space of functions $f$ locally in 
$\rL^2(\R^d)$ (resp. $\rL^2(\R^d\setminus \{0\})$) such that 
\begin{equation}\label{norm-B}
\| f\|_{\rB}\ := \ \| f\|_c\, +\, \sum_{j=1}^\infty 2^{j/2}\| f\|_{c_j}
\hspace{.4cm}
\Bigl(\mbox{resp.}\hspace{.2cm}\| f\|_{\dot{\rB}}\ := \ \sum_{j\in \Z}2^{j/2}\| f\|_{c_j}\Bigr)
\end{equation}
is finite. Its dual $\rB^\ast$ (resp. $\dot{\rB}^\ast$) is equipped with 
\begin{equation}\label{norm-B*}
\| f\|_{\rB^\ast}\ := \ \max\bigl(\| f\|_c\, ; \, \sup_{j\geq 1}2^{-j/2}\| f\|_{c_j}\bigr)
\hspace{.4cm}
\Bigl(\mbox{resp.}\hspace{.2cm}\| f\|_{\dot{\rB}^\ast}\ := \ \sup _{j\in \Z}2^{-j/2}\| f\|_{c_j}
\Bigr)\period 
\end{equation}
One can easily check that the embeddings $\rL _s^2\subset \rB\subset 
\rL _{1/2}^2$, for any $s>1/2$, and $\rB\subset \dot{\rB}$, are all continuous. 
Notice that, for $\rS=\rL _s^2, \rB$, and $\dot{\rB}$,
\begin{equation}\label{duality}
\forall f\in \rS^\ast, \, \forall g\in \rS,\,  \bar{f}g\in
\rL^1\hspace{.3cm}\mbox{ and }\hspace{.3cm}
\bigl|\langle f\, ,\, g\rangle \bigr|\ \leq\ \|f\|_{\rS^\ast}\cdot \|g\|_{\rS}\period 
\end{equation}
For $z\in \rho (P(h))$, $R(z;h)$ can be viewed as a bounded operator from $\rL _s^2$ 
to $\rL _{-s}^2$, for $s\geq 0$, and from $\rB$ to $\rB^\ast$, being a bounded operator on 
$\rL^2(\R^d)$. 
When $d\geq 3$, one can show using Hardy's inequality (\ref{hardy}) 
that, for $z\in \rho (P(h))$, $R(z;h)$ can even be viewed as a bounded operator from $\dot{\rB}$ to 
$\dot{\rB}^\ast$ (cf. \cite{wz}), a stronger result. \\
Let $I$ be a compact interval included in $]0;+\infty[$ and $d\geq 3$. By 
\cite{fh}, we know that $I$ contains no eigenvalue of $P(h)$. By Mourre's commutator theory (cf. \cite{abg,m}), we also know that for fixed $h$, $\|R(\cdot;h)\|_{\rS ,\rS^\ast}$ is 
bounded on $\{z\in \C; \Re z\in I, \Im z\neq 0\}$
whenever $\rS=\rL _s^2$ ($s>1/2$) or $\rS=\rB$. Adapting an argument by \cite{wz}, 
the above norm is even seen bounded when $\rS=\dot{\rB}$.  Summarizing, for $s>1/2$ and any given $h>0$,
the following chain of inequalities holds true
\begin{equation}\label{order}
\sup_{\stackrel{\Re z\in I}{\Im z\neq 0}}\, \|R(z;h)\|_{\rL _s^2,\rL _{-s}^2}\ \leq \ 
\sup_{\stackrel{\Re z\in I}{\Im z\neq 0}}\, \|R(z;h)\|_{\rB ,\rB^\ast}\ \leq \ 
\sup_{\stackrel{\Re z\in I}{\Im z\neq 0}}\, \|R(z;h)\|_{\dot{\rB} ,\dot{\rB}^\ast}\ < \ \infty
\period
\end{equation}
\subsection{The Non-Trapping Condition.}

We now estimate the terms involved in (\ref{order}) as $h \to 0$. When
$V=0$, it is known that
\begin{equation}\label{sup-resolv}
\sup_{\stackrel{\Re z\in I}{\Im z\neq 0}}\, \|R(z;h)\|_{\rS ,\rS^\ast}\ = \ O(1/h)\comma
\end{equation}
whenever $\rS=\rL _s^2$ ($s>1/2$), or $\rS=\rB$.
Our aim is to 
characterize those potential $V$ for which (\ref{sup-resolv}) holds
true with $\rS=\rB$. \\
If $V\in C^\infty(\R^d,\R)$
and satisfies (\ref{decroit}), then a characterization of those $V$'s such that (\ref{sup-resolv}) holds true
is well known, at least in the case
$\rS=\rL _s^2$ ($s>1/2$) or $\rS=\rB$, as we now describe. Let 
$T^\ast \R^d\ni (x,\xi)\donne p(x,\xi):=|\xi |^2+V(x)$ be the symbol 
of $P(h)$. Since the potential $V$ is bounded below, 
for all energies $\lambda$ on $p^{-1}(\lambda)$ the speed $|\xi|$ 
is bounded above. Thus the particle cannot escape to infinity in
finite time and
$p$ defines a complete smooth Hamilton flow
$(\phi^t)_{t\in \R}$ on $T^\ast \R^d$. The symbol $p$
is said {\em non-trapping} at the energy $\lambda$ whenever 
\begin{equation}\label{non-captif}
\forall (x,\xi) \in p^{-1}(\lambda)\comma \hspace{.5cm}\lim_{t\tend -\infty}|\phi^t(x,\xi)|
\, =\, +\infty \ \mbox{and}\ \lim_{t\tend +\infty}|\phi^t(x,\xi)|
\, =\, +\infty \period
\end{equation}
In many cases it is easy to show trapping by topological criteria, see
\cite{KK1}.\\
Let $\rS=\rL _s^2$ with $s>1/2$, or $\rS=\rB$. Then (\ref{sup-resolv}) holds true if 
and only if any energy $\lambda\in I$ is non-trapping for $p$
(cf. \cite{gm,jec2,ma2,rt,vz,w}). This statement 
has been extended to the homogeneous space $\rS=\dot{\rB}$ 
($d\geq 3$) by \cite{cj}, for $V$'s of class $C^2$ only. \\
First note that such a characterization is a Bohr correspondence principle: 
in the limit $h\to 0$, a qualitative property of the classical flow (the 
non-trapping condition) is connected to a propagation property of the quantum
evolution operator $U(t;h)=\exp (-ih^{-1}tP(h))$. 
Indeed the propagation estimate (\ref{int-propa}) turns out to be equivalent to the above estimate (\ref{sup-resolv}). \\
Second, it is also useful to develop a semiclassical, stationary 
scattering theory (the case $\rS=\rL _s^2$ actually suffices). If the non-trapping condition is true, one expects
to deduce from (\ref{sup-resolv}) 
bounds on several scattering objects (as is done when $V \in C^\infty(\R^d)$, cf. \cite{r2,rt}). If trapping occurs,
one expects that the resonant 
phenomena have a leading order influence on the scattering objects (cf. \cite{gm2,r2}).  \\
Of course, these two motivations are still present if 
Coulomb singularities are allowed.\\
When only repulsive Coulomb singularities occur, it was proved in \cite{kmw} that 
the non-trapping condition implies that  (\ref{sup-resolv}) is true with $\rS=\rL _s^2$ ($s>1/2$). If 
at least one attractive Coulomb singularity is present, the flow
is not complete anymore and the 
previous non-trapping condition does not even make sense. However, it is
known that one can "regularize the flow" 
(see \cite{kn} and references therein), and it turns out the
regularization is easier to deal with in 
dimension $d=2$ and $d=3$.\\
In the present paper, we choose to focus on the case $d=3$,
which is the physically important situation. Our study is devoted to
generalizing the previous characterization, in a case where the potential admits arbitrary Coulomb singularities. Note that
we do expect our results extend 
to the case $d=2$. \\
Let $d=3$ and assume that $\cS$ contains an attractive singularity. Let $(x,\xi )\in T^\ast\hat{M}=T^\ast(\R^3\setminus \cS)$. As we shall see 
in Subsection~\ref{regu-flow}, there exists some at most countable subset $\J(x,\xi ) 
\subset \R$ and a smooth function $\phi (\cdot ;x,\xi):\R\setminus \J(x,\xi )\dans 
T^\ast\hat{M}$ such that 
$\phi (\cdot ;x,\xi )$ solves the Hamilton equations generated by the symbol $p$ of $P(h)$ 
with initial value $(x,\xi)$ (see (\ref{cauchy-pb}) and (\ref{p-coulomb}) below). Furthermore, 
for all $t\in \R\setminus \J(x,\xi )$, $p(\phi (t;x,\xi))=p(x,\xi)$. The function 
$\phi$ replaces the usual flow.
It is thus natural to say that $p$ is non-trapping at 
energy $\lambda$ whenever
\begin{equation}\label{non-trapping}
\forall (x,\xi) \in p^{-1}(\lambda)\comma \hspace{.5cm}\lim_{t\tend
  -\infty}|\pi _x\phi (t;x,\xi)|
\, =\, +\infty \ \mbox{and}\ \lim_{t\tend +\infty}|\pi _x\phi (t;x,\xi)|
\, =\, +\infty \comma 
\end{equation}
where $\pi _x\phi (t;x,\xi)$ denotes the configuration or base component of $\phi (t;x,\xi)\in T^\ast\hat{M}$. 

\subsection{Survey.}
In view of (\ref{order}) and (\ref{non-trapping}), we can now state 
our main result. 
\begin{theorem}\label{main}
Let $V$ be a potential satisfying the assumptions (\ref{decroit}) and (\ref{sing-sj}). If 
there are no attractive singularities $(N=0)$, let $d\geq 3$ else let $d=3$. 
Let $I_0$ be an open interval included in $]0;+\infty[$. The following 
properties are equivalent. 
\begin{enumerate}
\item For all $\lambda \in I_0$, $p$ is non-trapping at energy $\lambda$. 
\item For any compact interval $I\subset I_0$, there exists $C>0$ such that, for $h\in ]0;h_\ast]$, 
\begin{equation}\label{size-resolvent}
\sup_{\stackrel{\Re z\in I}{\Im z\neq 0}}\, \|R(z;h)\|_{\rB ,\rB^\ast}
\ \leq \ C\, h^{-1}\period 
\end{equation}
\end{enumerate}
\end{theorem}
In \cite{w3}, the point $2$ of Theorem~\ref{main} is derived from a
virial-like assumption, which is stronger than the non-trapping
condition. It is assumed there that only one singularity occurs and
that (\ref{sing-sj}) holds true for a constant $f_j$. The statement ``$1\impl 2$'' 
of Theorem~\ref{main} is proved in \cite{kmw} when $N=0$. Theorem~\ref{main} provides the converse.
More importantly, it extends the result to the delicate case $N>0$. \\
To complete the picture given by Theorem~\ref{main}, we study in Section~\ref{when-non-trapping?} 
the non-trapping condition.
In the case of a {\em single} Coulomb singularity, we show that it is always
satisfied when the energy $\lambda$ is large enough, as in the case of a smooth potential
(see Remark \ref{rem:5.5}). The classically forbidden 
region in configuration space
then is a point (for attracting Coulomb potential), or it is diffeomorphic
to a ball (in the repelling case). 
Conversely Proposition \ref{trapping} says that -- irrespective 
of the number of singularities and the energy ---
{\em only} for the case of a single point or ball trapping does not
need to occur.
In particular, Corollary~\ref{effective-trapping} states that trapping always occurs
for two or more singularities at large enough energies.\\
We point out that our proof of Theorem~\ref{main} gives some additional insight about the case 
when the non-trapping condition fails at some energy $\lambda>0$ 
(cf. Proposition~\ref{almost-invariant}). In such a situation, ``semiclassical trapping'' 
occurs, as described by (\ref{B-bounds}) and (\ref{L2-bounds}). Notice that a resonance phenomenon 
(cf. \cite{hs,lm}) is a particular case of the quasi-resonance phenomenon defined in \cite{gs}, the 
latter being a particular case of our ``semiclassical trapping'' criterion. Propositions~\ref{supp-maj} 
and~\ref{almost-invariant} show that the ``semiclassical trapping'' is microlocalized near ``trapped trajectories" (see (\ref{traj-bornee}) for a precise definition). It would be interesting to check whether a 
(quasi-)resonance phenomenon is related to our ``semiclassical trapping'' 
(cf. Remark~\ref{resonant-state}). A traditional study of the resonances 
``created" by a bounded trajectory  (see \cite{gs} and references therein) would also be of interest.
We do hope that the present paper may help to overcome the difficulties due to the singularities. \\
While the proof of ``$2\impl 1$" in Theorem~\ref{main} follows the strategy developed by \cite{w} 
for smooth potentials, we use a rather different argument compared to \cite{rt,gm,kmw} when 
showing ``$1\impl 2$". In these papers, a semiclassical version of Mourre's commutator theory is used 
(cf. \cite{abg,m}), and the Besov-like space $B$ is replaced by the weaker $\rL _s^2$ ($s>1/2$). An 
alternative approach is given in \cite{b} for compactly supported perturbations
of the Laplacian, using a contradiction argument due to G. Lebeau in \cite{l}. This method was adapted 
in \cite{jec2} to include long-range, smooth perturbations, the study still being carried out in the space 
$\rL _s^2$ ($s>1/2$). This technique
was further developed in \cite{cj} to tackle the estimates in the optimal homogeneous space $\dot{\rB}$, 
by combining and adapting an original estimate derived in \cite{pv}. Note that both works \cite{pv} 
and \cite{cj} only require $C^1$ resp. $C^2$ smoothness on the potential. Note also that
the extension of Theorem \ref{main} to the homogeneous estimate in $\dot{\rB}$ still is open.
Now, the contradiction argument of \cite{jec2,cj}
is a key ingredient of the present study. Concerning the treatment of the
singularities, we stress that our study uses many 
results from \cite{gk}, the propagation results being here
crucial. The main features we need on the 
regularization of the classical flow are provided by \cite{gk,kn}. Our
main new contributions are given in Proposition~\ref{evolution-u_h} and in
Section~\ref{general-case}.  \\
Finally, we give some nonrelativistic, physical situations for which
our result applies. In both examples below, we may add to the operator a smooth exterior
potential satisfying (\ref{decroit}). 
\begin{example}\label{ex:mol0}
The behaviour of a particle with charge $e_0$ in the presence of fixed, point\-li\-ke ions,
with nonzero charges $z_1,\ldots,z_{N'}$, is governed by the operator (here $d=3$)
\begin{equation}
P_1(h)\ :=\ -h^2\Delta_x+\sum_{j=1}^{N'}\frac{e_0z_j}{|x-s_j|}. 
\label{Mol:op}
\end{equation}
The hydrogen atom corresponds to $N'=1$, $z_1>0$, and $e_0<0$. Clearly (\ref{decroit}) and
(\ref{sing-sj}) hold true. If charges have different sign, the
model has attractive and repulsive singularities. 
\end{example}
\begin{example}\label{ex:mol}
Consider a molecule with $N'$ nuclei having positive charges $z_1,\ldots,z_{N'}$, binding 
$K>0$ electrons with charge $-1$. We assume the nuclei are fixed
(Born-Oppenheimer idealization), and we neglect electron-electron
repulsion. The behaviour of each electron  is then  
governed by $P_1(h)$ in (\ref{Mol:op}). Let $h_0>0$ be fixed. 
Let $\psi_k$ be the normalized wavefunction of $P_1(h_0)$ of electron number $k$. 
Let $\rho _k=|\psi _k|^2$ be its charge density. Consider another,
much heavier particle with charge $e_0$. Its scattering by the
molecule can be described by $P(h)$ where 
\begin{equation}
V(x):=e_0\left(\sum_{j=1}^{N'}\frac{z_j}{|x-s_j|}+\sum_{k=1}^{K}W_k(x)\right),
\mbox{ with }
W_k(x):=-\int_{{\R}^d}\frac{\rho_k(q)}{|q-x|}\,dq\period  
\label{VW}
\end{equation}
As we show in Section \ref{sectonex}, it turns out that the $\psi_k$'s are ``nice enough'' to make $W_k$
well defined, smooth away from the singularities
$s_1,\ldots,s_{N'}$, and to make $W_k$ satisfy (\ref{decroit}). Though
(\ref{sing-sj}) does not hold, we show the proof of our result applies in 
this case. 
\end{example}
%

\section{Preliminaries.}
\label{pseudo-sing}
\setcounter{equation}{0}

We shall often use well known facts concerning $h$-pseudodifferential 
calculus, functional calculus, and semiclassical measures in the sequel. For sake of
completeness, we recall here the main 
results we need, referring to \cite{dg,g,gl,h,ke,lp,ma1,n,r1} for further details. Since our 
Schr\"odinger operator has Coulomb singularities, it does not define a pseudodifferential operator yet. 
For this reason, we also explain here how we can use pseudodifferential calculus ``away from the 
singularities'':
the required results are essentially contained in \cite{gk}; notice however that we do not need the results 
in the appendix of \cite{kmw}, which are, by the way, not known if an attractive Coulomb 
singularity is present. Last, we also recall basic results on the regularization of the Hamilton flow 
when an attractive singularity is present, refering to \cite{gk,ke,kn} for details. 

\subsection{Symbolic calculus with singularities.}
\label{symb-calculus}

Let $d\in \N^\ast$. For $(r,m)\in \R^2$, we consider the vector space (space of symbols)
\begin{eqnarray}\label{symb-class}
\Sigma_{r;m}&:=&\Big\{a\in C^\infty(T^*\R^d) \ ;\ 
\forall \gamma=(\gamma_x,\gamma_\xi)\in \N^{2d}\comma \, \exists C_\gamma>0\, ;\, \hfill \\
&&\hspace*{30mm}\sup_{(x,\xi )\in T^\ast\R^d}\, \langle x\rangle
^{-r+|\gamma_x|}\, \langle 
\xi\rangle 
^{-m+|\gamma_\xi|}\bigl| (\partial^\gamma a)(x,\xi )\bigr| \ \leq \ C_\gamma\Big\}\period 
\nonumber
\end{eqnarray}
If $r,m\leq 0$, then $\Sigma_{r;m}$ is contained in the
vector space of bounded symbols, which are 
smooth functions $a:T^\ast\R^d\dans \C$ such that
\begin{equation}\label{bd-symb-class}
\forall \gamma\in \N^{2d}\comma \, \exists C_\gamma>0\, ;\, 
\sup_{(x,\xi )\in T^\ast\R^d}\, \bigl| (\partial^\gamma a)(x,\xi )\bigr| \ \leq \ C_\gamma\period 
\end{equation}
For a larger class of symbols $a$, one can define the {\em Weyl $h$-quantization} of $a$, 
denoted by $a^w_h$. It acts on $u\in C_0^\infty(\R^d)$ as follows (cf. \cite{dg,ma1,n,r1}). 
\begin{equation}\label{w}
\bigl(a^w_hu\bigr)(x)\, = \,  (2\pi h)^{-d}\int_{\R^d} e^{i\xi\cdot (x-y)/h}a\bigl((x+y)/2,\xi\bigr)u(y)\, 
dy\,d\xi \period 
\end{equation}
If $a$ is a bounded symbol, then $a_h^w$ extends to a bounded operator on 
$\rL^2(\R^d)$, uniformly with respect to 
$h$, by Calder\'{o}n-Vaillancourt's theorem (cf.\ \cite{dg,ma1,r1}).
We shall also use the following functional calculus of Helffer-Sj\"ostrand, which can 
be found in \cite{dg,ma1}. Given $\theta\in C_0^\infty(\R)$, one can 
construct an almost analytic extension $\theta^{\C}\in C_0^\infty(\C)$
(with $\overline{\partial}\theta^{\C}(z)={\cal O}(\Im(z)^\infty))$. 
Let $H$ be a self-adjoint operator in some Hilbert space. The bounded operator $\theta (H)$, defined 
by the functional calculus of self-adjoint operators, can be written as 
\begin{equation}\label{hs}
\theta (H)\, = \,  \frac{-1}{\pi}\int_{\C} \frac{d\theta^{\C}}{d\barre{z}}(z)\cdot 
(z-H)^{-1}\, d\cL_2(z)\period 
\end{equation}
where $\cL_2$ denotes the Lebesgue mass on $\C$. \\
Let us now recall some well known facts about semiclassical measures, which can be found in \cite{g,gl,ke,lp}. 
Let $(u_n)_n$ be a bounded sequence in $\rL^2(\R^d)$. Up to extracting a subsequence, we may
assume that it is pure, i.e.\ it has a unique semiclassical measure $\mu $. By definition $\mu $ is a finite, nonnegative
Radon measure on the cotangent space $T^\ast\R^d$. Furthermore, there exists a sequence $h_n\tend 0$ such 
that, for any $a\in C_0^\infty(T^\ast\R^d)$, 
\begin{equation}\label{caract-measure}
\lim_{n\tend \infty}\, \bigl\langle u_n\, ,\, a^w_{h_n}u_n\bigr\rangle\ = \ 
\int_{T^\ast\R^d}a(x,\xi )\, \mu (dx\, d\xi)\, =:\, \mu (a)\period 
\end{equation}
One may relate the total mass of 
$\mu$ to the $\rL^2$-norm of the $u_n$'s (see \cite{gl}, or \cite{ke,lp}),
through the following
\begin{proposition}[\cite{gl}]\label{total-mass}
Let $(u_n)_n$ be a pure bounded sequence in $\rL^2(\R^d)$ such that 
\begin{eqnarray}
\lim_{R\to +\infty}\limsup_{n\to\infty}\int_{|x|\geq R}|u_n(x)|^2\, dx&=&0\comma \label{x-mass}\\
\lim_{R\to +\infty}\limsup_{n\to\infty}\int_{|\xi|\geq R/h_n}|\cF u_n\, (x)|^2\, d\xi&=&0
\comma\label{xi-mass}
\end{eqnarray}
where $\cF u_n$ denotes the Fourier transform of $u_n$. Then the sequence $(\|u_n\|^2)_n$ 
converges to the total mass $\mu (T ^\ast\R^d)$ of its semiclassical measure $\mu$. 
\end{proposition}
\Pf See the proof of Proposition 1.6 in \cite{gl}.\qed

Besides, transformation of the semiclassical measure upon composition 
of the $u_n$'s with a diffeomorphism is described in the
\begin{proposition}[\cite{gl}]\label{diffeo}
Let $\Phi :U\dans V$ be a $C^1$ diffeomorphism between two open subsets 
of $\R^p$ ($p\geq 1$). Let $\Phi _c :T^\ast U\dans T^\ast V$ be the symplectomorphism
\begin{equation}\label{canonical-map}
(y,\eta )\ \donne \ \bigl(\Phi (y)\, ;\, (\Phi '(y)^T)^{-1}\eta \bigr)\period 
\end{equation}
Here $\Phi '(y)^T$ denotes the transpose of $\Phi '(y)$. 
Given $a\in C_0^\infty(T^\ast V)$, let $b\in C_0^\infty(T^\ast U)$ 
be defined by $b=a\circ \Phi _c$. Then, for every compact subset $K$ of $V$, 
\[\lim_{h\to 0}\, \sup _{\stackrel{\|u\|\leq 1}{\support u\subset K}}\, \bigl\|
(a^w_hu)\circ \Phi  \, -\, b^w_h(u\circ \Phi )\bigr\|\ = \ 0\period\]
Let $K$ be a compact subset of $V$ and $(u_n)_n$ be a pure bounded sequence in $\rL^2(V)$ such 
that, for all $n$, $\support u_n\subset K$. Denote by $\mu$ its semiclassical 
measure. Then the sequence $(u_n\circ \Phi)_n$ is bounded in $\rL^2(U)$, its semiclassical measure $\tilde{\mu}$ is given by $|{\rm Det} \Phi '|^{-1}\Phi _c ^{-1}(\mu)$, and $\mu (a)=\tilde{\mu} (b)$.
\end{proposition}
\Pf See the proof of Lemma 1.10 in \cite{gl}.\qed

We now focus on the treatment of Coulomb singularities in dimension
$d\geq 3$, in combination with the $h$-pseudodif\-ferential framework. 
To begin with, let us recall Hardy's inequality.
\begin{equation}\label{hardy}
\forall f\in C_0^\infty(\R^d)\comma \ \int_{\R^d}\frac{|f(x)|^2}{|x|^2}\, dx\ 
 \leq \ \frac{4}{(d-2)^2}\, \|\nabla _xf\| ^2 \ =\ \frac{4}{h^2(d-2)^2}\, \| h\nabla _xf\|^2, 
\end{equation}
where the last bound is relevant in the present, semiclassical regime.\hfill\\
We next discuss how one can use  $h$-pseudodifferential calculus "away from the singularities". 
Recall that $\hat{M}=\R^d\setminus\cS$.
Let $\chi\in C_0^\infty(\R^d)$ with $\chi =1$ near the set $\cS$ of all singularities 
. Define the (truncated) $h$-pseudodifferential operator 
\begin{equation}\label{regularized-P}
P_\chi(h)\ :=\ -h^2\Delta _x+(1-\chi )V\period 
\end{equation}
Its symbol 
\begin{equation}\label{symb-regularized-P}
T^\ast\R^d\ni (x,\xi)\donne p_\chi (x,\xi)=|\xi |^2+\bigl(1-\chi (x)\bigr)V(x)
\end{equation}
belongs to $\Sigma _{0;2}$ (cf. (\ref{symb-class})). 
The following lemma is essentially proved in \cite{gk}. 
\begin{lemma}\label{pseudo-away-sing}
Let $d\geq 3$. Let $\chi\in C_0^\infty(\R^d)$ with $\chi =1$ near 
$\cS$ 
and $\theta\in C_0^\infty(\R)$. Let $P_\chi (h)$ be given by (\ref{regularized-P}).
Let $T>0$, $k,k'\in \R$, $r,m\in\R$, and $[-T;T]\ni t\donne a(t)\in \Sigma _{r;m}$ be a 
continuous function such that, for all $t\in [-T;T]$, $a(t)=0$ near
$\support \chi$.\\ 
Then, in $C^0\big([-T;T];\cL(\rL _k^2;\rL _{k'}^2)\big)$,  
\begin{eqnarray}
\bigl(P(h)-P_\chi (h)\bigr)\bigl(a(\cdot )\bigr)_h^w &=&O(h^2)\comma \label{P-P_chi}\\
\mbox{and, if }m\leq 2\comma \ \bigl(\theta (P(h))-\theta(P_\chi(h))\bigr)\bigl(a(\cdot )\bigr)_h^w &=&O(h^2)\period 
\label{theta-P-P_chi}
\end{eqnarray}
\end{lemma}
\Pf Let $r,m,k,k'\in \R$. For $a\in \Sigma _{r;m}$ and $f\in
\cS(\R^d)$, the Schwartz space on $\R^d$,   
\[\langle \cdot\rangle ^{k'}\bigl(P(h)-P_\chi (h)\bigr)a_h^w\langle \cdot \rangle ^kf
\ =\ V(-h^2\Delta+1)^{-1}\cdot (-h^2\Delta+1)\chi \langle \cdot\rangle ^{k'}a_h^w
\langle \cdot \rangle ^kf\comma\]
where $V(-h^2\Delta+1)^{-1}\in \cL(\rL^2;\rL^2)$ has norm 
$O(1/h^2)$ by (\ref{hardy}). Now, if $a$ is replaced by a continuous map $t\mapsto a(t)$ with
$a(t)=0$ near $\support \chi$ for all $t$, then, for all $N\in \N$, 
\[(-h^2\Delta +1)\chi \langle \cdot\rangle ^{k'}a(\ast )_h^w\langle \cdot \rangle ^k 
= O\bigl(h^N\bigr)\]
in $C^0([-T;T];\cL(\rL^2;\rL^2))$, by the usual $h$-pseudodifferential calculus. This yields 
(\ref{P-P_chi}). On the other hand
it is known that, for all $k\in \R$, the resolvents $(P(h)+i)^{-1}$ and 
$(P_\chi (h)+i)^{-1}$ are bounded from $\rL _k^2$ to $\rL _k^2$ 
(see \cite{rs4}, Sect. XIII.8), and,
by (\ref{hardy}), there exists some $\alpha (k)\geq 0$ such that 
\[\bigl\| (P(h)+i)^{-1}\bigr\| _{\cL (\rL _k^2;\rL _k^2)}\ = \ O\bigl(h^{-\alpha (k)}\bigr)
\comma 
\ \bigl\| (P_\chi (h)+i)^{-1}\bigr\| _{\cL (\rL _k^2;\rL _k^2)}\ = \ O\bigl(h^{0}\bigr)\period \]
Besides, there is a 
$\chi _1\in C^\infty(\R^d)$ with $\chi\chi_1=0$, $\chi _1=1$ at infinity, 
and $\chi _1a(t)=a(t)$ for all $t$.
Hence, for all $N\in \N$, $(1-\chi _1)a(\ast )_h^w\langle \cdot \rangle ^k=O(h^N)$ in 
$C^0\big([-T;T];\cL(\rL^2;\rL^2)\big)$. We may now adapt the arguments 
in the proof of Lemma 3.1 in \cite{gk} to get (\ref{theta-P-P_chi}). \qed

\subsection{Extension of the flow.}
\label{regu-flow}

Here we explain how the usual flow can be extended when attractive
singularities occur (more details are given in \cite{ke,kn}). \\
Let $d=3$. We still denote by $p$ the smooth function defined by 
\begin{equation}\label{p-coulomb}
p:\ \hat{P}\dans \R\qmbox{,}
(x,\xi )\donne |\xi |^2\, +\, V(x)\qmbox{where}\hat{P}:= T^\ast\hat{M} \period
\end{equation}
Let $\pi _x$ (resp. $\pi _\xi$) be the projection $T^\ast\R ^d\dans \R ^d$ 
defined by $\pi _x(x,\xi ):=x$ (resp. $\pi _\xi(x,\xi ):=\xi$).
As for any smooth dynamical system, the hamiltonian initial value problem,
\begin{align}\label{cauchy-pb}
\nonumber
&
\frac{dX}{dt}(t)=\nabla _\xi p\, (X(t);\Xi (t))\mbox{ , }\frac{d\Xi}{dt}(t)= 
-\nabla _x p\, (X(t);\Xi (t)),\\
&
(X(0);\Xi (0))\, =\, \ph =(x,\xi )\in \hat{P}
\end{align}
has a unique maximal solution $\phi:\hat{D}\ar \hat{P}$ with 
\[\hat{D} = \left\{(t,\ph)\in{\R}\times \hat{P}\, ;\, 
t\in\ ]T^-(\ph),T^+(\ph)[\ \right\},\]
where the functions $T^\pm:\hat{P}\ar \overline{{\R}}$ satisfy $T^- < 0 < T^+$
and are lower resp.\ upper semi-continuous with respect to the 
natural topology on the extended line $\overline{{\R}}:=\{-\infty\}\cup{\R}\cup\{+\infty\}$. 
In particular, the set $\hat{D}\subseteq{\R}\times \hat{P}$ is open.\\
If no attractive singularity is present (i.e.\ $N=0$ in the notation of Paragraph \ref{1.1}),
then $\hat{D} = {\R}\times \hat{P}$.
Otherwise a maximal solution can fall on an attractive 
singularity $s$ at finite time $T^+(\ph)>0$. 
Such a time is called a {\em collision time}. 
In that case, it turns out that,
setting 
\begin{eqnarray*}
\J(\ph)&:=&
\left\{\begin{array}{ccc}
\emptyset&\mbox{if}&T^-(\ph)=-\infty,\ T^+(\ph)=\infty\\
\{T^+(\ph)\}&\mbox{if}&T^-(\ph)=-\infty,\ T^+(\ph)<\infty\\
\{T^-(\ph)\}&\mbox{if}&T^-(\ph)>-\infty,\ T^+(\ph)=\infty\\
\{T^+(\ph)\}+{\mathbb Z}\Big(T^+(\ph)-T^-(\ph)\Big)&\mbox{if}&
T^-(\ph)>-\infty,\ T^+(\ph)<\infty
\end{array}\right. \\
\mbox{and}\hspace{.2cm} D&:=& \left\{(t,\ph)\in{\R}\times \hat{P}\, ;\, 
t\not \in\ \J(\ph)\ \right\}\comma 
\end{eqnarray*}
the map $\phi$ can be uniquely extended to a smooth map $D\ar
\hat{P}$, still denoted by $\phi$.
Even more, when $T^+(\ph)<\infty$, {\em backscattering} occurs, that
is, for $0<t<T^+(\ph)-T^-(\ph)$, we have
\begin{eqnarray}
\pi _x\phi(T^+(\ph)+t;\ph) &=&\pi _x\phi(T^+(\ph)-t;\ph), \nonumber \\
\pi_\xi\phi(T^+(\ph)+t;\ph) &=& -\pi_\xi\phi(T^+(\ph)-t;\ph),\label{time-reversal}
\end{eqnarray}
and one may set $\pi _x\phi (T^+(\ph);\ph )=s$. We mention that the momentum
$\pi _\xi\phi (\cdot ;\ph )$ however blows up at $T^+(\ph)$, in the
following sense:
\[\lim _{t\to T_+(\ph )}|\pi _\xi\phi (t;\ph )|=\infty\comma\mbox{ while }
v :=\lim _{t\nearrow T_+(\ph )}\frac{\pi _\xi\phi (t;\ph )}{|\pi _\xi\phi
(t;\ph )|}=-\lim _{t\searrow T_+(\ph )}\frac{\pi _\xi\phi (t;\ph )}{|\pi _\xi\phi
(t;\ph )|}\mbox{ exists}\period \] 
For any $\ph \in \hat{P}$, we obtain in
this way a configuration trajectory $(\pi _x\phi (t;\ph ))_{t\in \R}$, 
which has a countable set 
$\J(\ph )$ of collision times $t_0$ for which 
\begin{equation}\label{collision-time}
\lim _{t\to t_0}\pi _x\phi (t;\ph )\in \{s_j, 1\leq j\leq N\}\hspace{.5cm}
\mbox{and}\hspace{.5cm}
\lim _{t\to t_0}|\pi _\xi\phi (t;\ph )|=\infty\period
\end{equation}
Although $\phi$ is not a complete flow on $\hat{P}$, 
the {\em broken trajectory} $(\phi (t;\ph ))_{t\in \R\setminus \J(\ph )}$ 
is a solution of (\ref{cauchy-pb}) on 
$\R\setminus \J(\ph )$. Its values lie in the energy shell $p^{-1}(p(\ph ))$. Note that
no collision with the repulsive singularities can occur. \\
For $t\in\R$, it is convenient to introduce $\phi ^t:D_t\ar \hat{P}$ defined by 
\begin{equation}\label{phi^t}
D_t\ :=\  \left\{\ph\in\hat{P}\, ;\, 
t\not \in\ \J (\ph)\ \right\}\hspace{.5cm}
\mbox{and}\hspace{.5cm}\phi ^t(\ph )\ :=\ \phi (t;\ph )\period 
\end{equation}
Note further that the Hamiltonian system $(\hat{P}, \omega_0, p)$ with canonical
symplectic form $\omega_0$ can be uniquely extended to a smooth Hamiltonian system with a 
complete flow (see Section \ref{when-non-trapping?}). \\
An important feature to analyse the pseudo-flow $\phi$ is the Kustaanheimo-Stiefel 
transformation (KS-transform for short). We briefly describe it here and refer to 
\cite{gk,ke,kn,ss}, for further details. For $z=(z_0,z_1,z_2,z_3)^T\in \R^4$, let 
\[\Lambda (z)\ =\ \left(
\begin{array}{cccc}
z_0&-z_1&-z_2&z_3\\
z_1&z_0&-z_3&-z_2\\
z_2&z_3&z_0&z_1
\end{array}
\right)\period \]
Let $\cK:\R^4\dans \R^3$ be defined by 
\begin{equation}\label{def-KS}
\cK (z)\ :=\ \Lambda (z)\cdot z
=\bsm z_0^2-z_1^2-z_2^2+z_3^3\\ 2z_0z_1-2z_2z_3\\2z_0z_2+2z_1z_3\esm
\period\hspace{.5cm}\mbox{For all }z\in \R^4\mbox{ it satisfies }\
|\cK (z)|=|z|^2\period 
\end{equation}
We call it the {\em Hopf map}.
See the Appendix for more information.\\
Let $\R ^3_\pm :=\{(x_1;x_2;x_3)\in \R^3; \pm x_1>0\}$ and $z\in \R^4$. It turns out (see \cite{gk}) that, if $x:=\cK (z)\in \R ^3_+$, $\cA _+(z):=\sqrt{2}(x_1+|z|^2)^{-1/2}(z_1+iz_4)\in S^1$ and, if $x:=\cK (z)\in \R ^3_-$, $\cA _-(z):=\sqrt{2}(-x_1+|z|^2)^{-1/2}(z_2+iz_3)\in S^1$. Furthermore, one can explicitly construct smooth maps $\cJ _\pm:\R ^3_\pm\times S^1\dans \R^4$ such that, locally,  
\begin{equation}\label{local-diffeo}
(\cK ,\cA _\pm)\circ \cJ_\pm\ =\ {\rm Id}\hspace{.2cm}\mbox{in}\hspace{.2cm}\R^3_\pm\times S^1\hspace{.4cm}
\mbox{and}\hspace{.4cm}\cJ _\pm\circ (\cK ,\cA _\pm)\ =\ {\rm Id}\hspace{.2cm}\mbox{in}\hspace{.2cm}\cJ
_\pm(\R ^3_\pm\times S^1)\period 
\end{equation}
For $z=\cJ _\pm (x;\theta )$, for $x\in \R^3_\pm$ and $\theta \in S^1$, 
we have $dz=C|x|^{-1}dx\,d\theta$ for some constant $C>0$. In
particular, there exists $C'>0$ such that, for all $f,g:\R^3\dans \C$
measurable, 
\begin{equation}\label{change-KS}
\int_{\R^3} |x|^{-1}\cdot |f(x)g(x)|\, dx\ =\ C'\, \int_{\R^4}
|f\circ \cK (z)\, g\circ \cK(z)|\, dz\period 
\end{equation}
It is useful to consider the following extension to 
phase space. For $\phz =(z;\zeta )\in T^\ast\R ^4$, we set as usual
$\pi_z\phz=z$ and $\pi_\zeta\phz=\zeta$. 
If $(x;\xi)\in T^\ast (\R ^3\setminus\{0\})$, let $z\in \R^4$ such that 
$x=\cK(z)=\Lambda (z)\cdot z$ ($z$ is not unique). Then, we define
\begin{equation}\label{sol-zeta}
\zeta \ := \ 2 \Lambda (z)^T \xi = 
2\left(\begin{array}{ccc}
z_0&z_1&z_2\\
-z_1&z_0&z_3\\
-z_2&-z_3&z_0\\
z_3&-z_2&z_1
\end{array}\right)
\left(\begin{array}{c}
\xi _1\\
\xi _2\\
\xi _3
\end{array}\right), 
\end{equation}
which is a solution of the equation $2|x|\xi =\Lambda (z)\zeta$. 
The {\em KS-transform} is defined by 
\begin{equation}\label{Lambda-ast}
\cK ^\ast:T^\ast (\R ^4\setminus\{0\})\dans T^\ast (\R
^3\setminus\{0\})\mbox{ , }\ \cK ^\ast (z;\zeta )\ = \ \Bigl(\Lambda (z)\cdot z\, ;\,
\frac{1}{2|z|^2}\Lambda (z)\cdot \zeta\Bigr)\period
\end{equation}
%

Assume that an attractive singularity sits at $0$. 
Recall that, by (\ref{sing-sj}), $V(x)=f(x)/|x|+W(x)$ on $\Omega\setminus\{0\}$ , where 
$\Omega :=\{x\in \R^3; |x|<r\}$ for some $r>0$, with $f,W\in C_0^\infty(\R^d ;\R)$. 
Let $\tilde\Omega :=\cK^{-1}(\Omega)$. 
Let $\ph_0=(x_0;\xi _0)\in \hat{P}$ be such 
that the first collision of $(\pi _x\phi (t;\ph _0))_{t\in \R}$ takes place at $0$ at time $t_+(\ph _0)>0$. 
Let $\cT _0$ be the connected component of $\big\{t\in \R;\ \pi _x\phi (t;\ph _0)\in \Omega\big\}$ containing $t_+(\ph _0)$. 
Let $z_0\in \R^4$ be such that $x_0=\cK(z_0)$ and let $\zeta _0$ be the $\zeta$ given by (\ref{sol-zeta}) with $(z;\xi )=(z_0;\xi _0)$. 
For $\pht =(t;\tau )\in T^\ast\R$, $\phz =(z;\zeta )\in T^\ast\R ^4$, 
let $\tilde{p}(\pht ;\phz ):=|\zeta |^2+f\circ \cK (z)+|z|^2(W\circ \cK (z)-\tau )$. 
Since $\tilde{p}$ is smooth on $T^\ast\R\times T^\ast\R ^4=T^\ast (\R_t\times\R _z^4)$, independent 
of $t$, and since its Hamilton vector field at point $(t,\tau ; z, \zeta)$ is given by $\l(-|z|^2,
0;2\zeta ,2\tau z\ri)$ outside a compact region in $(z, \zeta)$, 
there exists a unique maximal solution $\R\ni s\donne \big(t(s);\tau (s);z(s);\zeta (s)\big)=
(\pht (s);\phz (s))$ to the Hamilton equations associated with $\tilde{p}$ 
\begin{equation}\label{cauchy-pb-extended}
\left(
\begin{array}{rclcrcl}
(dz/ds)(s)& =& \nabla _\zeta \tilde{p}\, (\pht (s);\phz (s))&\hspace{.5cm}
&(d\zeta /ds)(s)& = & -\nabla _z \tilde{p}\, (\pht (s);\phz (s))\\
(dt/ds)(s)& =& \nabla _\tau \tilde{p}\, (\pht (s);\phz (s))&\hspace{.5cm}
&(d\tau /ds)(s)& = & -\nabla _t \tilde{p}\, (\pht (s);\phz (s))
\end{array}
\right)\comma 
\end{equation}
with initial condition $(\pht (0);\phz (0))=(\pht _1;\phz _1)$. We denote it by 
\[\tilde{\phi}(s;\pht _1;\phz _1)\ :=\ \bigl(\pht (s;\pht _1;\phz _1);\phz (s;\pht _1;\phz _1)\bigr)
\ =\ \bigl(t(s;\pht _1;\phz _1);\tau (s;\pht _1;\phz _1);z(s;\pht _1;\phz _1);
\zeta (s;\pht _1;\phz _1)\bigr) \period  \]
Let $(\pht _0;\phz _0)=(0;p(\ph _0);z_0;\zeta _0)$. It turns out that, 
for all $t_1\in\cT _0$, there exists a unique $s\in \R$ such that 
$t_1=t(s;\pht _0;\phz _0)$. Furthermore, if $t_1\neq t_+(\ph _0)$, $z(s;\pht _0;\phz _0)\neq 0$ and  
\begin{equation}\label{phi-phitilde}
\phi (t_1;\ph _0)\ =\ \phi \bigl(t(s;\pht _0;\phz _0);\ph _0\bigr)\ =\ 
\cK ^\ast \bigl(\phz (s;\pht _0;\phz _0)\bigr)\period 
\end{equation}
%

\section{Towards the non-trapping condition.}
\label{resolv-esti-non-trapp}
\setcounter{equation}{0}

The aim of this section is to prove the implication ``$2\impl 1$'' of Theorem~\ref{main}. 
We thus assume that $2$ holds true and we want to show that $p$ is non-trapping 
at energy $\lambda$, for all $\lambda \in I_0$. Let $\lambda _0$ be such an energy.
We can find a compact interval $I\subset I_0$ such that $\lambda _0$
belongs to the interior of $I$. 
By assumption, (\ref{size-resolvent}) holds true for $I$. This 
implies, by (\ref{order}), that (\ref{sup-resolv}) holds true for $\rS=\rL _s^2$, for any 
$s>1/2$. As in \cite{jec3}, we follow the strategy in \cite{w}. We 
translate the bound on the resolvent into a bound on a time integral of the associated propagator 
\begin{equation}\label{propa}
U(t;h)\ :=\ \exp (-ih^{-1}tP(h))\period 
\end{equation}
If an attractive singularity is present (and $d=3$), we need some information on the time-dependent 
microlocalization of $U(t;h)u_h$, for some family $(u_h)_h$ of $\rL^2(\R^3)$ functions. 
Most of it is already available in \cite{gk}. We also need some well-known facts on the classical 
flow, which we borrow from \cite{kn}. In Subsection~\ref{co-st-evo}, we shall recall results from
\cite{gk,kn} and extend them a little bit. Then we proceed with the announced proof in 
Subsection~\ref{nece-non-trap}. In the repulsive case ($N=0$ and $d\geq 3$), we show in 
Subsection~\ref{repulsive-case} that Wang's proof may be carried over with minor changes.

\subsection{Coherent states evolution.}
\label{co-st-evo}

In this subsection we are interested in the case where an attractive
singularity occurs (i.e.\ $N>0$) but the results hold true for
$N=0$. Better results in the latter case are given in
Subsection~\ref{repulsive-case}. Proposition~\ref{evolution-u_h} is 
the main result of the subsection.

Before considering the time evolution of coherent states, we recall some 
basic facts on the classical dynamics, in particular on the dilation function
\[a_0:T ^\ast\R^d\dans \R\ \mbox{ , }\ a_0(x,\xi ):=x\cdot \xi\period \]
\begin{lemma}\label{escape-infinity}
Consider a dimension $d\geq 2$ and energies $\lambda>0$.
\begin{enumerate}
\item
Then for some $R_1=R_1(\lambda)\ge R_0$ and all 
$\ph _0:=(x_0,\xi_0)\in p^{-1}(]\lambda/2 ;\infty [)$
\begin{eqnarray}\label{imply}
|x_0|\geq R_1\ \impl \  \{p,a_0\} (\ph _0)\geq \lambda
/2\ \comma \mbox{ and}\\
\label{result-kn}
\liminf_{t\to\pm\infty}
 | \pi _x\phi (t;\ph _0)|>
R_1\ \impl \  \lim_{t\to \pm\infty}|\pi _x\phi
(t;\ph_0)|=+\infty\period
\end{eqnarray}
\item
For any $T,R>0$, there is some $R_2>R_1$ such that, 
for all $\ph _0=(x_0,\xi_0)\in p^{-1}(]\lambda/2 ;2\lambda[)$
\begin{equation}\label{loca-far-away}
 \bigl(|x_0|>R_2\bigr)\ \impl \ 
\bigl(|  \pi _x\phi (t;\ph _0)| >R\mbox{ for all }t\in [-T;T]\bigr)\period 
\end{equation}
\end{enumerate}
\end{lemma}
\Pf
We shortly recall the standard arguments. 
Thanks to the decay properties (\ref{decroit}) of $V$,
\[\{p,a_0\}(\ph _0)= 2\big(p(\ph _0)-V(x_0)\big)-\langle x_0,\nabla
V(x_0)\rangle \ge \lambda/2\]
for large $|x_0|$, implying (\ref{imply}). As the dilation function $a_0$
is the time derivative of the phase space function $|x|^2/2$, composed
with $\phi$, the second time derivative of the
latter function is eventually bounded below by $\lambda/2>0$, 
if the l.h.s.\ of (\ref{result-kn}) is satisfied. Thus
$t\mapsto |\pi _x\phi(t;\ph_0)|^2$ goes to infinity, showing (\ref{result-kn}). 
Let $V_0=\inf_{|x|\ge R_0} V(x)$. 
Relation (\ref{loca-far-away}) follows, since the speed is bounded
by $|\xi_0|\le (4\lambda-2V_0)^{1/2}<\infty$.
\qed

For $h\in ]0;h_\ast]$ the {\em dilation operator} $E_h$ on 
$\rL^2(\R^d)$, given by
$E_h(f)(x):=h^{-d/4}f(h^{-1/2}x)$, is unitary, as are the 
{\em Weyl operators} 
\[w(\ph_0 ;h)\ :=\ \exp \bigl(ih^{-1/2}(x_0\cdot x-\xi_0\cdot D_x)\bigr)\qquad\mbox{ for }
\ph _0:=(x_0,\xi_0)\in T^*\R^d, \]
(cf.\ \cite{h}, p.~151, \cite{fo}). 
The {\em coherent states operators}, 
microlocalized at $\ph_0$, are
\begin{equation}\label{op-etat-co}
c(\ph_0 ;h)\ :=\ E_h \cdot w(\ph_0 ;h)
\end{equation}
A direct computation gives that
\begin{eqnarray*}
E_h ^\ast a^w_h E_h &=& \bigl( a(h^{1/2}\star ;h^{-1/2}\star)\bigr)_1^w\comma\\
c(\ph_0 ;h)^\ast a^w_hc(\ph_0 ;h) &= &\bigl( a(x_0+h^{1/2}\star ;
\xi _0+h^{-1/2}\star)\bigr)_1^w\comma 
\end{eqnarray*}
where $b(\star ;\star)$ denotes the symbol $(x;\xi)\mapsto b(x;\xi )$. 
It is known (cf. \cite{w0}) that 
\begin{equation}\label{symb-coherent-state}
\forall a\in \Sigma _{0,0}\comma\forall f\in \cS(\R^d)\comma \, c(\ph_0 ;h)^\ast a^w_hc(\ph_0 ;h)
f\ =\ a(\ph_0)f\, +\, O(h)\comma
\end{equation}
where $\cS (\R^d)$ denotes the Schwartz space on $\R^d$. 
Let $u_h$ be the function given by 
\begin{equation}\label{translated-gaussian}
u_h\ :=\ c(\ph_0 ;h)\, \pi^{-d/4}\, \exp \bigl(-|\cdot |^2/2\bigr)\period 
\end{equation}
Then $(u_h)_h$ is a family of $\rL^2(\R^d)$-normalized coherent states microlocalized 
at $\ph _0$. We collect properties of the family $(U(\cdot ;h)u_h)_h$
of the propagated states. \\
In the remainder part of Subsection~\ref{co-st-evo} we consider initial conditions
in phase space 
\[\ph _0:=(x_0,\xi_0)\in\hat{P} \qmbox{with energy}\lambda:= p(\ph _0)>0\] 
and the associated coherent states $(u_h)_h$ microlocalized 
at $\ph _0$.\\
In \cite{gk} the following {\em energy localization} of $(u_h)_h$ is obtained. 
We give a short proof using Lemma~\ref{pseudo-away-sing}.
\begin{lemma}[\cite{gk}]\label{u_h-energy-localization}
Let $d\geq 3$ and $\theta\in C_0^\infty(\R)$ such that $\theta =1$
near $\lambda$. Then, in $\rL^2(\R^d)$, 
$(1-\theta (P(h)))u_h=O(h)$. 
\end{lemma}
\Pf Let $\chi ,\tilde{\chi}\in C_0^\infty(\R^d)$ with $\chi ,\tilde{\chi}=1$ near
$\cS$, $\chi ,\tilde{\chi} =0$ near $x_0=\pi_x\ph _0$, and 
$\chi \tilde{\chi}=\chi$. From (\ref{translated-gaussian}), we see that $\tilde{\chi}u_h=O(h)$ 
in $\rL^2(\R^d)$. By Lemma~\ref{pseudo-away-sing}, 
\[\bigl(1-\theta (P(h))\bigr)u_h\ =\ \bigl(1-\theta (P(h))\bigr)(1-\tilde{\chi})u_h \, +\, O(h) 
\ =\ \bigl(1-\theta (P_\chi (h))\bigr)(1-\tilde{\chi})u_h \, +\, O(h)\comma \]
in $\rL^2(\R^d)$, where $P_\chi (h)$ is as in
(\ref{regularized-P}). Besides, thanks to (\ref{symb-coherent-state}) and
using (\ref{symb-regularized-P}), 
\begin{eqnarray*}
\bigl(1-\theta (P(h))\bigr)u_h &=& \bigl(1-\theta (P_\chi (h))\bigr)
u_h\, +\, O(h) \ =\ \bigl(1-\theta (p_\chi (\ph _0))\bigr) u_h\, +\,
O(h)\\
&=&0\, +\,  O(h) \period \qed
\end{eqnarray*}
 From \cite{gk} we pick the following {\em localization away from singularities}
\begin{lemma}[\cite{gk}]\label{sing-evolution-u_h}
Let $d=3$. Let $K$ be a compact subset of $\R$ such that $K\cap \J(\ph _0)=\emptyset$ 
(cf. (\ref{collision-time})). If $\sigma \in C^\infty_0(\R^3)$ has small enough a support near 
the set $\cS$ of singularities, then $K\ni t\donne \sigma U(t;h)u_h$ 
is of order $O(h)$ in $C^0\big(K;\rL^2(\R^3)\big)$. 
\end{lemma}
\Pf See the proof of Theorem 1, p.\ 25 in \cite{gk}.\qed

A careful inspection of the result in \cite{gk} on the frequency set
shows that even after a collision, we have the following 
{\em localization along our broken trajectories}:
\begin{lemma}\label{time-space-localization}
Let $d=3$. 
Let $K$ be a compact subset of $\R$ such that $K\cap \J(\ph _0)=\emptyset$ 
(cf. (\ref{collision-time})). 
Let $\epsilon >0$ and $K\ni t\donne a(t;\ast )\in C^\infty_0
(\hat{P})$ be continuous functions such that 
$a(t;x,\xi )=0$ if $|x-\pi _x\phi (t;\ph _0)|\leq \epsilon$. 
Then $(a(\cdot;\ast))_h^wU(\cdot ;h)u_h=O(h)$ in $C^0\big(K;\rL^2(\R^3)\big)$.
\end{lemma}
\Pf See the proof of Theorem 1, p. 25 in \cite{gk}.\qed

We also need to complete Lemma~\ref{time-space-localization} with a {\em bound} on 
$(U(\cdot ;h)u_h)_h$ {\em near infinity in position space} and, since the singularities are far away, we can assume $d\geq 3$. This is the purpose of the following 
\begin{lemma}\label{egorov}
Let $d\geq 3$. Let $T>0$ and $R:=\max (R_0;1+|x_0|\}$. Let $R_2>R_1$ large enough 
such that (\ref{loca-far-away}) holds true. Let $R_3>R_2+1$ and 
$\kappa\in C^\infty(\R^d;\R)$ such that $\support\kappa\subset 
\{y\in \R^d;\,  |y|>R_2+1\}$ and $\kappa =1$ on $\{y\in \R^d;\, |y|>R_3\}$. Then 
\[\kappa U(\cdot ;h)u_h \ =\ O(h)\qmbox{in} C^0\big([-T;T];\rL^2(\R^3)\big).\] 
\end{lemma}
\Pf 
The proof is based on an Egorov type estimate which is valid
although $P(h)$ is not a pseudodifferential operator.\\
$\bullet$ Let $\tau\in C_0^\infty(\R^d)$ such that $\tau =1$ on $\{y\in \R^d; |y|\leq R_0\}$ and $\tau =0$ near the set 
\[\pi _x\, \bigcup _{t\in [-T;T]}\, \left(p^{-1}\bigl(]\lambda /2  ;2\lambda [\bigr)\cap
\phi ^t\bigl(\{(x,\xi);\, |x|>R_2\}\bigr)\right)\period \]
This is well-defined by (\ref{phi^t}), (\ref{loca-far-away}), and the choice
of $R$. Let $p_\tau$ be defined as in (\ref{symb-regularized-P}).
Let $\theta\in C_0^\infty(\R)$ with $\support \theta \subset ]\lambda /2 
;2\lambda [$ such that $\theta =1$ near $\lambda$. Set 
\begin{equation}\label{def:a}
  a:T ^\ast\R^d\dans \C\qmbox{,}
  a(x,\xi )=\kappa(x) \, \theta(p_\tau (x,\xi)).
\end{equation}  
Thanks to (\ref{loca-far-away}), 
$[-T;T]\ni t \donne a\circ \phi ^t$ is a $\Sigma _{0;0}$-valued, $C^1$-function. Therefore, by Calder\'{o}n-Vaillancourt $(a\circ \phi ^{t})_h^w$ is $h$--uniformly bounded, and for 
$t\in [-T;T]$, strongly in $\rH^2(\R^d)$, 
\begin{eqnarray}
\lefteqn{U(t;h)^\ast a_h^w U(t;h)\, -\, \bigl(a\circ \phi ^t\bigr)_h^w
\quad=\quad\int_0^t\, \frac{d}{ds}\Bigl( U(s;h)^\ast \bigl(a\circ \phi ^{t-s}\bigr)_h^wU(s;h)
\Bigr)\, ds} \label{ansatz-egorov}\\
&=& \int_0^t\, U(s;h)^\ast\ \left(\frac{i}{h}\bigl[P(h),
\bigl(a\circ \phi ^{t-s}\bigr)_h^w\bigr]\, +\, 
\bigl((d/ds)a\circ \phi ^{t-s}\bigr)_h^w\right)\ U(s;h)\, ds \period\nonumber
\end{eqnarray}
The support properties of $a$ and the choice of $\tau$ ensure that, for all $r\in [-T;T]$, 
$(d/dr)a\circ \phi ^{r}=\{p,a\circ \phi ^{r}\}=\{p_\tau ,a\circ 
\phi ^{r}\}$. Thus, (\ref{ansatz-egorov}) equals 
\[\int _0^t\, U(s;h)^\ast\ 
\Bigl(\frac{i}{h}\bigl[P(h),\bigl(a\circ \phi ^{t-s}\bigr)_h^w\bigr]\, -\, 
\bigl(\{p_\tau ,a\circ \phi ^{t-s}\}\bigr)_h^w\Bigr)\ U(s;h)\, ds\period \]
By Lemma~\ref{pseudo-away-sing}, (\ref{ansatz-egorov}) equals 
\[\int _0^t\, U(s;h)^\ast\ 
\Bigl(\frac{i}{h}\bigl[P_\tau(h),\bigl(a\circ \phi ^{t-s}\bigr)_h^w\bigr]\, -\, 
\bigl(\{p_\tau ,a\circ \phi ^{t-s}\}\bigr)_h^w\Bigr)\ U(s;h)\, ds\, +\, hB_h(t) \]
where $[-T;T]\ni t\donne B_h(t)\in \cL(\rL^2(\R^d))$ is bounded, 
uniformly with respect to $h$. 
By the usual pseudodifferential calculus, 
\[[-T;T]\ni r\ \donne \ \frac{i}{h}\bigl[P_\tau(h),\bigl(a\circ \phi ^{r}\bigr)_h^w\bigr]\, -\, 
\bigl(\{p_\tau ,a\circ \phi ^{r}\}\bigr)_h^w\in \cL(\rL^2(\R^d))\]
is $O(h)$ in $C^0\big([-T;T];\cL(\rL^2(\R^d))\big)$. 
Thus, so is (\ref{ansatz-egorov}).\hfill\\
$\bullet$ 
Since  $a\circ \phi ^{t}$ vanishes near $x_0$, for $t\in [-T;T]$, 
$t\donne (a\circ \phi ^t)_h^w u_h$ is $O(h)$ in $C^0\big([-T;T];L^2\big)$, 
by (\ref{symb-coherent-state}). 
Thus so is $ U(\cdot ;h) \left(a\circ\phi^t\right)_h^w u_h$. \\
$\bullet$ Using the previous points, the Lemmata~\ref{u_h-energy-localization} and~\ref{pseudo-away-sing}, 
the fact that $\theta (P(h))$ and $U(t,h)$ commute, and the usual pseudodifferential calculus,
\begin{eqnarray}
\hspace{-0.4cm}
\kappa U(t ;h)u_h&=&
\kappa U(t ;h)\theta (P(h))u_h \, +\, O(h)\ =\ \kappa \theta(P(h)) U(t ;h)u_h \, +\, O(h)\nonumber\\
&=&\kappa \theta(P_\tau(h)) U(t ;h) u_h \, +\, O(h)
\ =\ (\kappa \theta(p_\tau))^w_h \, U(t ;h) u_h \, +\, O(h)
\nonumber\\
&=& U(t ;h) \, \left(\kappa \theta(p_\tau)\circ \phi^t \right)^w_h  u_h \, +\, O(h)\ =\ O(h)\period \qed
\label{kappa-energy-localization}
\end{eqnarray}

 From these lemmata, we can deduce the following information on the time evolution of the 
coherent states $u_h$. 
\begin{proposition}\label{evolution-u_h}
Let $N>0$ and $d=3$.
Let $K$ be a compact subset of $\R$ such that $K\cap \J(\ph _0)=\emptyset$ 
(cf. (\ref{collision-time})). Let $\tau\in C^\infty_0(\R^3)$ with $\tau =1$ near $0$. 
For $t\in \R$ and $x\in \R^3$, set $\tau_t(x) :=
\tau (x-\pi _x\phi (t;\ph _0))$. Take the support of $\tau$ small 
enough such that, for all $t\in K$, 
$\support(\tau _t) \cap \cS=\emptyset$. 
Then, for any $a\in \Sigma _{0;0}$ and any $t\in K$, 
\begin{eqnarray*}
a^w_h \, U(t;h) u_h
&=& (\tau _t a)^w_h \, U(t;h) u_h \, +\, e(t)\comma 
\end{eqnarray*}
where $K\ni t\donne e(t)$ is  $O(h)$ in $C^0(K;\rL^2(\R^d))$. 
\end{proposition}
\Pf Let $T>0$ such that $K\subset [-T;T]$. Let $\kappa_0,\kappa_1\in C^\infty(\R^3;\R)$ 
such that $\kappa_0+\kappa_1=1$ and $\kappa:=\kappa_1$ satisfies the assumptions 
of Lemma~\ref{egorov}. Then, by Lemma~\ref{egorov}, 
\begin{eqnarray*}
a^w_h \, U(t ;h) u_h &=& 
a^w_h \, \kappa_0 \, U(t ;h)u_h \, +\, O(h)\comma 
\end{eqnarray*}
in $C^0(K):=C^0(K;\rL^2(\R^d))$. Now let $\sigma_0\in C_0^\infty(\R^3;\R)$ such that 
$\sigma_0=1$ near $s_j$ for any $1\leq j\leq N'$, and, for all 
$t\in K$, $\support\sigma_0\cap\support\tau_t=\emptyset$. 
Upon possibly decreasing the support of $\sigma_0$, we may apply 
Lemma~\ref{sing-evolution-u_h}. This yields 
\begin{eqnarray*}
a^w_h \, U(t ;h) u_h
&=&
a^w_h \, \kappa_0 \, (1-\sigma_0) \,
U(t ;h)u_h \, +\, O(h)\comma
\end{eqnarray*}
in $C^0(K)$. Let 
$\sigma\in C_0^\infty(\R^3;\R)$ such that 
$\sigma=1$ near each singularity $s_j$ and $\sigma \sigma_0=\sigma$.
For an energy cutoff  $\theta$ as in
Lemma~\ref{u_h-energy-localization}, we obtain, 
as in the proof of Lemma~\ref{egorov}
(see (\ref{kappa-energy-localization})),
\begin{eqnarray*}
&&
a^w_h \, \kappa_0 \, (1-\sigma_0) \,
U(t;h)u_h =
a^w_h \, \kappa_0 \, (1-\sigma_0) \, \theta(P_{\sigma}(h)) \, U(t ;h)u_h
\, +\, O(h)\comma
\end{eqnarray*}
in $C^0(K)$, since $1-\sigma_0$ is localized away from the singularities. 
By pseudodifferential calculus, 
\[a^w_h \, U(t ;h)u_h\ =\
b^w_h \, U(t ;h) u_h \, +\, O(h)\comma\]
in $C^0(K)$, where $b:=\theta(p_{\sigma}) \, (1-\sigma_0) \, \kappa_0 \, a\in 
C^\infty_0(\hat{P})$. 
Applying Lemma~\ref{time-space-localization} to $a(t)=(1-\tau _t)b$, 
\begin{eqnarray*}
a^w_h \, U(t ;h)u_h
&=&  (\tau _{t}b)_h^w
U(t ;h)u_h \, +\, O(h)\ =\
\bigl((1-\sigma_0) \, \kappa_0 \, \tau _{t}a \bigr)_h^w
U(t ;h)u_h \, +\, O(h)\comma
\end{eqnarray*}
in $C^0(K)$. Since $\tau _t \, (1-\sigma_0) \, \kappa_0=\tau _t$, for all $t\in K$, 
we obtain the desired result.
\qed

\subsection{Necessity of the non-trapping condition.}
\label{nece-non-trap}

Assuming $N>0$ and $d=3$, we want to show that (\ref{size-resolvent}) 
implies the non-trapping condition, yielding the proof of ``$2\impl 1$''. The proof below 
actually works if $N=0$, but a more straightforward and easier proof is provided in 
Subsection~\ref{repulsive-case}.\\ 
In view of (\ref{order}), we assume (\ref{sup-resolv}) for $\rS=\rL _s^2$ with 
$s>1/2$. This means that, for any $\theta\in C_0^\infty(I_0;\R)$, 
$\langle \cdot\rangle ^{-s}\theta (P(h))$ is 
Kato smooth with respect to $P(h)$ (by Theorem XIII.30 in \cite{rs4}).  This can be formulated in the 
following way (cf.\ Theorem XIII.25 in \cite{rs4}). There exists $C_s>0$ such that for any 
$\theta\in C_0^\infty(I_0;\R)$, 
\begin{equation}\label{int-propa}
\forall u\in \rL^2(\R^d)\comma \ \int_\R \, \|\langle \cdot \rangle ^{-s}\, U(t;h)\, 
\theta (P(h))u\|^2\, dt\ \leq \ C_s\cdot \|u\|^2\period 
\end{equation}
uniformly in $h\in ]0;h_\ast]$. Take $\lambda\in I_0$ and a function $\theta\in C_0^\infty(I_0;\R)$ 
such that $\theta =1$ near $\lambda $. Let $\ph _0:=(x_0,\xi_0)\in p^{-1}(\lambda )$ and consider 
the coherent states $u_h$ given by (\ref{translated-gaussian}). Let $(t_j)_{j\in J}$, 
with $J\subset \N^\ast$, be the set of collision times of the broken trajectory 
$(\phi (t;\ph _0))_{t\in \R}$ (cf. (\ref{collision-time})). Eq. (\ref{int-propa}) implies that, for all 
$T>0$ and all $h\in ]0;h_\ast]$, 
\[\int_{[-T;T]} \, \|\langle \cdot\rangle ^{-s}\, U(t;h)\, \theta (P(h))u_h\|^2\, dt\ \leq \ C_s
\period\]
We know from Subsection~\ref{regu-flow} that the collision 
times in $\J(\ph _0)$ have positive minimal distance $T^+(\ph
_0)-T^-(\ph _0)$.
Thus we can choose $\epsilon>0$ smaller than one fourth of that
distance, and define for $T>0$ the compact sets
\[K(T):=\{t\in[-T,T]; {\rm dist}(t,\J(\ph _0))\ge \epsilon\} \period\]
Notice that the length of $K(T)$ goes to infinity when $T\to \infty$,
while
\[\int_{{K(T)}} \, 
\|\langle \cdot\rangle ^{-s}\, U(t;h)\, \theta (P(h))u_h\|^2\, dt\ \leq \ C_s
\period\]
By energy localization of the coherent state 
(Lemma~\ref{u_h-energy-localization}) and Pythagoras' theorem, 
\begin{equation}\label{bound-without-theta}
\int_{{K(T)}} \, \|\langle \cdot\rangle ^{-s}\, U(t;h)\, u_h\|^2\, dt\, +\, O_T(h)\ \leq \ 2C_s
\comma
\end{equation}
where $O_T(h)$ is a $T$-dependent $O(h)$. 
We apply Proposition~\ref{evolution-u_h} for the bounded symbol $(x,\xi)\donne a(x,\xi )=
\langle x\rangle ^{-2s}$ and the compact ${K(T)}$ introduced above. This yields 
\[\int_{{K(T)}} \, 
\bigl\langle U(t;h)u_h\, ,\, \tau _t \langle \cdot\rangle ^{-2s}U(t;h)u_h\bigr\rangle\, dt
\, +\, O_T(h)\ \leq \ 2C_s\period\]
We can require that the support of the function $\tau$ is so small that, for all $t\in {K(T)}$, 
$\support(\tau _t) \cap \cS=\emptyset$ and $\tau _t \langle \cdot
\rangle ^{-2s}\geq (1/2)\tau _t\langle \pi _x\phi (t;\ph _0)\rangle ^{-2s}$. 
Therefore, 
\[\int_{{K(T)}} \, 
\langle \pi _x\phi (t;\ph _0)\rangle ^{-2s}\, \bigl\langle U(t;h)u_h\, ,\, 
\tau _t U(t;h)u_h\bigr\rangle\, dt\, +\, O_T(h)\ \leq \ 4C_s\period\]
Now, we apply Proposition~\ref{evolution-u_h} again for the bounded symbol 
$(x,\xi)\donne a(x,\xi)=1$, yielding
\begin{equation}\label{int-flow+O(h)}
\int_{{K(T)}} \, 
\langle \pi _x\phi (t;\ph _0)\rangle ^{-2s}\, dt\, +\, O_T(h)
\ \leq \ 4C_s\comma
\end{equation}
since the $u_h$ are normalized. Letting $h$ tend to $0$, we obtain, for all $T>0$,  
\begin{equation}\label{int-flow}
\int_{{K(T)}} \, \langle \pi _x\phi (t;\ph _0)\rangle ^{-2s}\, dt\ \leq \ 4C_s\period
\end{equation}
Assume semi-boundedness of the trajectory, that is, 
for some $t_0\in \R$, 
\begin{equation}\label{trapped-traj}
\{\pi _x\phi (t;\ph _0)\comma \, \pm t\geq t_0\}\ \subset\ 
\{y\in\R^3; |y|\leq R_1 \}\comma
\end{equation}
then, by (\ref{int-flow}), $4C_s$ is larger than $R_1^{-2s}$ times the length of
\[{K(T)} \setminus \{t\in \R ; \pm t< t_0\}\ = \ [-T;T]\setminus 
\Bigl\{t\in \R; \pm t< t_0\mbox{ and } {\rm dist}(t,\J(\ph _0))<\epsilon\Bigr\} \period\]
This is a contradiction since the latter tends to $\infty$ as $T\to \infty$. Thus (\ref{trapped-traj}) 
is false and we can apply (\ref{result-kn}), yielding 
the non-trapping condition (\ref{non-trapping}).

\subsection{The repulsive case.}
\label{repulsive-case}

Here we consider the case where any singularity is repulsive (i.e.\ $N=0$) and $d\geq 3$.  
We want to show that (\ref{size-resolvent}) implies the non-trapping condition.
Thanks to Proposition~\ref{egorov-N=0} below, we show that Wang's proof can be followed 
in the present case, yielding a much simpler proof than the one in
Subsection~\ref{nece-non-trap}. 

First of all, we show that an important ingredient in Wang's proof is 
available, namely the following weak version of Egorov's theorem. 
\begin{proposition}\label{egorov-N=0}
Let $N=0$ and $d\geq 3$. Let $T>0$ and $a\in \Sigma _{0;0}$. 
Let $\theta ,\gamma\in C^\infty_0(\R)$ such that $\gamma \theta =\theta$. 
Then $[-T;T]\ni t \donne \gamma (p)(a\circ \phi ^t)$ is a $\Sigma _{0;0}$-valued, $C^1$-function. 
Furthermore, there exists $C>0$, depending on $\theta$ and $a$, such that, for any $\epsilon >0$, for any $t\in [-T;T]$, 
\[U(t;h)^\ast a^w_hU(t;h)\theta (P(h)) \ = \ \bigl((\gamma (p)(a\circ \phi ^t))^w_h\, +\, r(t)\bigr)\theta (P(h))\comma \]
where $[-T;T]\ni t\donne r(t)$ is bounded by $C\epsilon +O_{\epsilon ,T}(h)$ in $C^0\big([-T;T];\cL(\rL^2(\R^d))\big)$. 
\end{proposition}
\Pf Let $\epsilon >0$. Since the singularities are repulsive, there exists some 
$\sigma_0\in C^\infty_0(\R^d,[0,1])$ which equals $1$ near each singularity, 
such that, $\epsilon ^2V\geq 1$ on the support of $\sigma_0$ and 
$(\sigma_0\circ\pi_x)(\gamma\circ p)=0$. 
Thus, for $g\in \rL^2(\R^d)$ and $f=\theta (P(h))U(t;h)g$, 
\beqno
\lefteqn{\|\sigma_0f\|^2\ \leq \ \langle f\, ,\, \sigma_0^2\epsilon ^2Vf\rangle \, +\, \epsilon
^2\langle \sigma_0f\, ,\, -h^2\Delta \, \sigma_0f\rangle}\\
&&\leq \ \epsilon ^2\langle \sigma_0^2f\, ,\, P(h)f\rangle+\epsilon ^2\langle \sigma_0f\, ,\, [-h^2\Delta ,\sigma_0]f\rangle\ \leq \ C_\theta ^2\epsilon ^2\|f\|^2+\epsilon ^2\langle \sigma_0f\, ,\, [-h^2\Delta ,\sigma_0]f\rangle \comma 
\eeqno
where $C_\theta$ depends only on $\theta$. 
Since $[-h^2\Delta ,\sigma_0]\theta (P(h))=O_\epsilon(h)$ in $\cL(\rL^2(\R^d))$,
\begin{equation}\label{near-sing}
\|\sigma_0\theta (P(h))U(t;h)\|\ \leq \ C_\theta\epsilon \, +\, O_\epsilon (h)
\end{equation}
in $C^0\big([-T;T];\cL(\rL^2(\R^d))\big)$. Let $\sigma\in C_0^\infty(\R^d)$ with $\sigma =1$ 
near each singularity such that $\sigma \sigma_0=\sigma$. Using (\ref{near-sing}), Lemma~\ref{pseudo-away-sing}, and pseudodifferential calculus, 
\begin{eqnarray*}
U(t;h)^\ast a^w_hU(t;h)\theta (P(h))&=&U(t;h)^\ast \bigl(a(1-\sigma_0)\bigr)^w_h\gamma (P_\sigma (h))\, U(t;h)\theta (P(h))\, +\, r_1(t)\\
&=&U(t;h)^\ast \bigl(a(1-\sigma_0)\gamma (p)\bigr)^w_h\, U(t;h)\theta (P(h))\, +\, r_2(t)\comma 
\end{eqnarray*}
where the $r_j$ are bounded by $C\epsilon +O_{\epsilon }(h)$ in $C^0\big([-T;T];\cL(\rL^2(\R^d))\big)$. By the choice of $\sigma_0$, $a(1-\sigma_0)\gamma (p)=a\gamma (p)=:a_\gamma$. Furthermore, for all $t\in [-T;T]$, $a_\gamma\circ \phi ^{t}=\gamma (p)(a\circ \phi ^{t})$ and $(d/dt)a_\gamma\circ \phi ^{t}=\{p,a_\gamma\circ \phi ^{t}\}=\{p_\sigma ,a_\gamma\circ \phi ^{t}\}$. This allows us to follow the arguments
in the proof of Lemma~\ref{egorov} showing that (\ref{ansatz-egorov}) with $a=a_\gamma$ is $O_{\epsilon
,T}(h)$ in $C^0\big([-T;T];\cL(\rL^2(\R^d))\big)$. \qed

Let $\lambda \in I_0$. As in Subsection~\ref{nece-non-trap}, (\ref{size-resolvent}) implies the existence of some constant $C_s>0$ such that (\ref{int-propa}) holds true, for $\theta\in C_0^\infty(I_0;\R)$ with $\theta (\lambda )=1$. Since no collision occurs, 
we choose $K(T)=[-T;T]$, take $a:(x,\xi)\donne \langle x\rangle^{-2s}$, and write (\ref{bound-without-theta}) as 
\begin{equation}\label{bound-with-a}
\int_{[-T;T]} \, \bigl\langle U(t;h)u_h\, ,\,
a_h^wU(t;h)u_h\bigr\rangle\, dt\, +\, O_T(h)\ \leq \ 2C_s
\period 
\end{equation}
By Lemma~\ref{u_h-energy-localization}, Proposition~\ref{egorov-N=0} with $\epsilon =C_s$, and (\ref{symb-coherent-state}), 
\begin{eqnarray*}
\bigl\langle U(t;h)u_h\, ,\, a_h^wU(t;h)u_h\bigr\rangle &=&\langle u_h\, ,\, 
\bigl(\gamma (p)(a\circ \phi ^t)\bigr)^w_hu_h\rangle\, +\, b_1(t)\ =\ \langle \pi _x\phi (t;\ph _0)\rangle ^{-2s}\, +\, b_2(t)
\end{eqnarray*}
where the $b_j$ are bounded by $CC_s+O(h)$ in $C^0([-T,T])$. This yields (\ref{int-flow+O(h)}), with bound $4C_s$ replaced by $(2+C)C_s$, and the non-trapping condition as in Subsection~\ref{nece-non-trap}.

\section{Semiclassical trapping.}
\label{semicl-trapp}
\setcounter{equation}{0}  

This section is devoted to the proof of the implication ``$1\impl 2$'' of Theorem~\ref{main}. 
We assume the non-trapping condition true on $I_0$ and we want to prove the bound 
(\ref{size-resolvent}), for any compact interval $I\subset I_0$. Here we follow the strategy in 
\cite{b,jec2}. We assume that the bound (\ref{size-resolvent}) is false, 
for some $I$. This means precisely that the following situation occurs, 
which we call ``semiclassical trapping''. There exist a sequence $(f_n)_n$ of 
nonzero functions of $\rH ^2(\R^d)$, a sequence 
$(h_n)_n\in ]0;h_0]^{\N}$ tending to zero, and a sequence $(z_n)_n\in \C^{\N}$ 
with $\Re(z_n)\tend \lambda \in I$ and $\Im(z_n)/h_n\to r\geq 0$, such that 
\begin{equation}\label{B-bounds}
\| f_n\|_{\rB^\ast}\ =\ 1\hspace{1cm}\mbox{and}\hspace{1cm}
\bigl\|(P(h_n)-z_n)f_n\bigr\|_{\rB}\ =\ o(h_n)\period 
\end{equation}
As in \cite{cj}, we shall see that "the $(f_n)_n$ has no $\rB^\ast$-mass at 
infinity" (see Proposition~\ref{compact-local} below). This yields the existence of some 
large $R_1'>0$, of a sequence $(g_n)_n$ of nonzero functions of $\rH ^2(\R^d)$, of a 
sequence $(h_n)_n\in ]0;h_0]^{\N}$ tending to zero, and of a sequence $(\lambda_n)_n
\in \R^{\N}$ with $\lambda_n\tend \lambda\in I$, such that 
\begin{equation}\label{L2-bounds}
\support g_n\subset \{x\in \R^3; |x|\leq R_1'\}\comma \hspace{.5cm}\| g_n\|\ =\ 1\comma 
\hspace{.5cm}\mbox{and}\hspace{.5cm}
\bigl\|(P(h_n)-\lambda_n)g_n\bigr\|\ =\ o(h_n)\period 
\end{equation}
Possibly after extraction of a subsequence, we may assume that the sequence $(g_n)_n$ has 
a unique semi-classical measure $\mu$, satisfying (\ref{caract-measure}) with $u_n$ replaced 
by $g_n$ (see Lemma \ref{g_n}).\\
Now, we look for a contradiction with the non-trapping condition. 
While, in the regular case, it is quite easy to show the invariance of $\mu$ under the 
flow generated by $p$, this is not clear in the present situation. We shall show the 
invariance for repulsive singularities in 
Subsection~\ref{repuls-sing}. In Subsection~\ref{general-case} however, we only show a weaker 
form of invariance, if there is an attractive singularity. This
Subsection~\ref{general-case} contains the main novelty of the paper. \\ 
The other steps of the strategy 
are essentially the same as in \cite{jec2}, 
as explained in Subsection~\ref{main-lines}. If the reader is only interested in the 
bound (\ref{size-resolvent}) with $B$ replaced by some $\rL _{s}^2$ ($s>1/2$),
we propose a simpler proof in Subsection~\ref{simpler-case}.

\subsection{Main lines of the proof.}
\label{main-lines}

In this subsection, we give the main steps leading to the contradiction between the "semiclassical trapping" and the non-trapping condition. Here we focus on the steps which are essentially proved as in \cite{jec2}. 
\begin{lemma}\label{r=0}
The sequence $\big(\|f_n\|^2\Im(z_n)/h_n\big)_n$ goes to $0$ and $\displaystyle\lim_{n\to\infty}\Im(z_n)/h_n=~0$.
\end{lemma}
\Pf We write $\|f_n\|^2\Im(z_n)=\Im \,\langle f_n\, ,\, (P(h_n)-z_n)f_n\rangle$, which is
$o(h_n)$ by (\ref{duality}) and (\ref{B-bounds}). This gives the first
result. Now, assume that $r>0$. Since $\|f_n\|^2(\Im(z_n)/h_n)$ goes to $0$,
$\|f_n\|$ must go to $0$, while 
$\|f_n\|\geq \| f_n\|_{\rB^\ast}=1$. This is a contradiction. \qed

Using (\ref{decroit}), we can show as in \cite{cj} the following localization in position space
\begin{proposition}\label{compact-local}
There exists $R_0'>R_0$ such that $\lim_{n\tend \infty}\, \|\un_{\{|\cdot|>R_0'\}}
f_n\|_{\rB^\ast}\ = \ 0$.
\end{proposition}
\Pf Let $a\in \Sigma _{0;0}$.
It is known that $(a_h^w)_{h\in ]0;h_\ast]}$ is uniformly bounded in $\cL (L^2_s;L^2_s)$
for any $s\in{\mathbb R}$. Even more,  
using a partition of unity adapted to the decomposition
${\mathbb R}^d=c\cup\left(\mathop{\cup}_{j\geq 1} c_j\right)$ from (\ref{c:c}), say
$1=\tau(x)+\sum_{j\geq 1} \tau_j(x)$, and writing, for any $u \in B^\ast$,
the identity $u = \tau \, u + \sum_{j \geq 1} \tau_j \, u$,
standard pseudodifferential calculus and almost orthogonality properties
allow to easily establish that 
$(a_h^w)_{h\in ]0;h_\ast]}$ is uniformly bounded in $\cL (B^\ast;B^\ast)$
(see \cite{cj} for a complete proof).
Now,
let $\alpha _n:=\langle f_n, ih_n^{-1}[P(h_n),a_{h_n}^w]f_n\rangle$. Expanding 
the commutator, using (\ref{duality}), (\ref{B-bounds}) and Lemma~\ref{r=0}, we observe that $\alpha _n\to 0$. 
For any $s>1/2$, $(f_n)_n$ is bounded in $\rL^2_{-s}$, 
since $\rL^2_{s}\subset B$. Now, we assume that $a$ vanishes near the set 
$\cS$ 
of all singularities. We can find $\chi\in C_0^\infty(\R ^d;\R)$ such that $a\chi =0$ and 
$\chi=1$ near the singularities. By Lemma~\ref{pseudo-away-sing} with $k=k'=-s$, 
\[\alpha _n \ = \ \langle f_n\, , \, ih_n^{-1}[P_\chi(h_n),a_{h_n}^w]f_n\rangle\, +\, O(h_n)
\period \]
Let $\theta \in C_0^\infty(\R;\R)$ with $\theta =1$ near $I$ and $\tilde\theta :=1-\theta$. 
Since $z_n\to \lambda\in I$, $(\|\tilde\theta (P(h_n))(P(h_n)-z_n)^{-1}\|)_n$ is uniformly 
bounded. Thus there exists $C>0$ such that 
\begin{eqnarray}
\|\tilde\theta (P(h_n))f_n\|_{B^\ast}&\leq & \max\bigl(\| \tilde\theta
(P(h_n))f_n\|_c\, ; \, \sup_{j\geq 1}2^{-j/2}\| \tilde\theta
(P(h_n))f_n\|_{c_j}\bigr)
\nonumber\\
&\leq &
C\|(P(h_n)-z_n)f_n\| = o(h_n)\comma \label{tilde-theta-Bast}
\end{eqnarray}
since (\ref{B-bounds}) implies that $\|(P(h_n)-z_n)f_n\|=o(h_n)$.
Using further that, for $s\in ]1/2;1]$, 
$\langle \cdot  \rangle ^{s}ih_n^{-1}[P_\chi(h_n),a_{h_n}^w]$ is uniformly bounded, 
\[\alpha _n \ = \ \langle f_n\, , \, ih_n^{-1}[P_\chi(h_n),a_{h_n}^w]\theta (P(h_n))f_n\rangle\, 
+\, O(h_n)\period \]
Since $ih_n^{-1}[P_\chi(h_n),a_{h_n}^w]$ is a $h$-pseudodifferential operator, we may 
apply Lemma~\ref{pseudo-away-sing} with $k=k'=-s$, yielding 
\[\alpha _n \ = \ \langle f_n\, , \, ih_n^{-1}[P_\chi(h_n),a_{h_n}^w]\theta (P_\chi(h_n))
f_n\rangle\, +\, O(h_n)\period \]
Using similar arguments again, we arrive at 
\begin{equation}\label{regular-terms}
\alpha _n \ = \ \langle \theta (P_\chi(h_n))f_n\, , \, ih_n^{-1}[P_\chi(h_n),a_{h_n}^w]
\theta (P_\chi(h_n))f_n\rangle\, +\, O(h_n)\period 
\end{equation}
Now we specify the symbol $a$ more carefully. By \cite{cj} (see
Proposition 8 and the second step of the proof of Proposition 7
therein), we can find $c>0$ and a function 
$\chi_1\in C_0^\infty(\R ^d)$ such that, for all $\beta =(\beta _j)\in \ell ^1$ with 
$|\beta |_{\ell ^1}=1$, there exists a symbol $a\in\Sigma _{0;0}$ satisfying the following 
properties. The function $\chi _1=1$ on a large enough neighbourhood of $0$ and of 
the support of $\chi $. The semi-norms of $a$ in $\Sigma _{0;0}$ are bounded independently of 
$\beta$ and, uniformly with respect to $\beta$,
\[\alpha _n \ \geq \ c\cdot \Bigl| \sum _j\beta _j 2^{-j}\bigl\|(1-\chi _1)\theta (P_\chi (h_n))f_n
\bigr\|_{c_j}^2\Bigr|\, +\, o(1)\period \]
By the above arguments, $\alpha _n\to 0$, uniformly in $\beta$. This implies that 
\[\sup _{j}\, 2^{-j}\bigl\|(1-\chi _1)\theta (P_\chi (h_n))f_n
\bigr\|_{c_j}^2\hspace{.3cm}\mbox{and therefore}\hspace{.3cm}
\sup _{j}\, 2^{-j/2}\bigl\|(1-\chi _1)\theta (P_\chi (h_n))f_n
\bigr\|_{c_j}\]
tend to $0$. In other words, $\|(1-\chi _1)\theta (P_\chi (h_n))f_n\|_{B^\ast}\to 0$. 
Since $B\subset \rL^2_{1/2-\epsilon}$ continuously, for any $\epsilon >0$, 
we derive from Lemma~\ref{pseudo-away-sing} that 
\[\bigl\|(1-\chi _1)\bigl(\theta (P_\chi (h_n))\, -\, \theta (P(h_n))\bigr)f_n
\bigr\|_{B^\ast}\ \to \ 0\comma \]
yielding $\|(1-\chi _1)f_n\|_{B^\ast}\to 0$, thanks to (\ref{tilde-theta-Bast}). 
Now the desired result follows for $R_0'$ large enough such that $|x|\geq R_0'\impl 
\chi _1(x)=0$.  \qed
\begin{lemma}\label{g_n}
Let $R_1'>R_0'$. There exist a sequence $(g_n)_n$ of nonzero functions of $\rH ^2(\R^d)$, 
bounded in $\rL^2(\R^d)$ and having a unique semiclassical measure $\mu$, a 
sequence $(h_n)_n\in ]0;h_\ast]^{\N}$ tending to zero, and a sequence 
$(\lambda_n)_n\in \R^{\N}$ with $\lambda_n\tend \lambda\in I$, 
such that (\ref{L2-bounds}) holds true. 
\end{lemma}
\Pf Let $\tau , \kappa\in C_0^\infty(\R ^d;\R)$ be such that 
$\support \tau ,\support \kappa\subset \{x\in \R ^d; |x|\leq R_1'\}$, 
$\kappa =1$ on $\{x\in \R^d; |x|\leq R_0'\}$, and $\tau \kappa 
=\kappa$. The sequence $(\tau f_n)_n$ is bounded in $\rL^2(\R^d)$. Possibly after 
extraction of a subsequence, we may assume that it has a unique semiclassical measure $\mu$. 
We shall show that 
\begin{align}
&\support \mu \subset \{(x,\xi)\in T^\ast\R ^d ; |x|\leq R_0'\}\comma \label{x-support}\\
& \support \mu \cap T^\ast(\R ^d\setminus \cS) \subset
p^{-1}(\lambda)\period
\label{energy-support}
\end{align}
By Proposition~\ref{compact-local}, $\|\un_{\{|\cdot|>R_0'\}}
\tau f_n\|$ goes to $0$. Using (\ref{caract-measure}), this implies
(\ref{x-support}). Now let $a\in C_0^\infty(T^\ast\R ^d)$ be such that
$a=0$ near $p^{-1}(\lambda)\cup \cS$. Since $(\|\langle \cdot\rangle^{-1}f_n\|)_n$ is bounded by 
(\ref{B-bounds}),
\begin{eqnarray}
\langle \tau f_n\, ,\,  a_{h_n}^w\tau f_n\rangle
&=&
\langle \tau f_n\, ,\,(\tau  a)_{h_n}^wf_n\rangle \, +\, O(h_n) 
\nonumber
\\
&=& \langle \tau f_n\, ,\, (\tau a)_{h_n}^w\theta (P(h_n))f_n\rangle \, 
+\, O(h_n)
\nonumber
\\
&& + \langle \tau f_n\, ,\, (\tau 
a)_{h_n}^w\tilde{\theta}(P(h_n))(P(h_n)-z_n)^{-1}\,(P(h_n)-z
_n)f_n\rangle \comma \label{energy-cut}
\end{eqnarray}
where $\theta \in C_0^\infty(\R;\R)$ with $\theta =1$ near $\lambda$, such that 
$\theta(p)a=0$, and $\tilde{\theta}=1-\theta$. By (\ref{B-bounds}), 
$\|(P(h_n)-z_n)f_n\|=o(h_n)$ and the last term in (\ref{energy-cut}) is a $o(h_n)$. 
We can find $\chi\in C_0^\infty(\R ^d;\R)$ such that $a\chi =0$ and $\chi=1$ near the singularities. By 
Lemma~\ref{pseudo-away-sing}, we recover
\[\langle \tau f_n\, ,\, a_{h_n}^w\tau f_n\rangle \ =\ \langle \tau f_n\, ,\, (\tau a)_{h_n}^w
\theta (P_\chi(h_n))f_n\rangle \, +\, O(h_n)\ =\ O(h_n)\]
since $a\theta (p_\chi)=0$. By (\ref{caract-measure}), this yields (\ref{energy-support}).\\
The symbol of $[-h_n^2\Delta ,\kappa ]$ belongs to $\Sigma_{-\infty,1}$ and 
is supported in $\{(x,\xi)\in T^\ast\R ^d ; R_0'<|x|<R'_1\}$. Let $\tilde\tau \in C^\infty_0(\R ^d)$
such that $\tilde\tau=1$ on $\support \nabla \kappa$ and $\support\tilde\tau\subset\{x\in \R ^d ; R_0'<|x|< R'_1\}$. Then $[-h_n^2\Delta ,\kappa ]f_n$
\begin{eqnarray*}
 &=&
[-h_n^2\Delta ,\kappa ]\tilde\tau f_n
\ =\
[-h_n^2\Delta ,\kappa ] \, (P_\chi (h_n)+i)^{-1} \,
(P_\chi (h_n)+i)\, \tilde\tau f_n
\\
&=&
[-h_n^2\Delta ,\kappa ] \, (P_\chi (h_n)+i)^{-1} \,
[-h_n^2\Delta ,\tilde\tau ]  f_n +
[-h_n^2\Delta ,\kappa ] \, (P_\chi (h_n)+i)^{-1} \,
\tilde\tau
\, (P(h_n)-z_n) \,  f_n
\\
&&
+
[-h_n^2\Delta ,\kappa ] \, (P_\chi (h_n)+i)^{-1} \, (i+z_n)
\, \tilde\tau f_n\ =:\ r_1+r_2+r_3 \period 
\end{eqnarray*}
Standard pseudodifferential calculus together with Proposition~\ref{compact-local} provide
$r_1=o(h_n^2)$, $r_2=o(h_n^2)$, and $r_3=o(h_n)$ in $\rL^2(\R^d)$.
Thus, setting $g_n:=\kappa f_n$,
\[(P(h_n)-z_n)g_n\ = \ \kappa (P(h_n)-z_n)f_n \, +\, o(h_n)
\ = \ o(h_n)\]
in $\rL^2(\R^d)$. By Proposition~\ref{compact-local} and (\ref{B-bounds}),
$\|g_n\|\to c$, with $c>0$, and $\Im(z_n)g_n=o(h_n)$ in
$\rL^2(\R^d)$, by
Lemma~\ref{r=0}. Setting $\lambda_n:=\Re (z_n)$, we obtain
$\|(P(h_n)-\lambda_n)g_n\|=o(h_n)$. Using (\ref{caract-measure}) and the previous arguments, $\mu$
is the unique semiclassical measure of $(g_n)_n$. \qed

We now collect properties of the $g_n$
and their semiclassical measure $\mu$, defined in Lemma \ref{g_n}.
\begin{lemma}\label{basic-properties}
Let $a\in C_0^\infty(T^\ast\R ^d)$ such that $a=0$ near the set $\cS$
of all singularities.\vspace*{-6mm}

\begin{enumerate}
\item
Then $\mu (\{p,a\})=0$ (``$\mu$ is invariant under the flow'').
\item
If $a=0$ near $p^{-1}(\lambda)$ or near $\{(x,\xi)\in T^\ast\R^d; |x|\leq R_0'\}$
then $\mu (a)=0$.
\item
Let $\tau \in C^\infty(\R ^d)$ such that $\tau =0$ near $\cS$. 
Then the sequence $(\|\tau ih_n\nabla g_n\|)_n$ is bounded. 
\end{enumerate}
\end{lemma}
\Pf 
{\bf 1)}
Let $\chi\in C_0^\infty(\R ^d;\R)$ such that $a\chi =0$ and $\chi=1$
near the set $\cS$ of all singularities.
In particular, $\{p,a\}=\{p_\chi ,a\}$. By Lemma~\ref{pseudo-away-sing}, 
\begin{eqnarray}\label{comm-expectation}
a_n\ :=\ \bigl\langle g_n\, ,\, ih_n^{-1}[P(h_n),a_{h_n}^w]g_n\bigr\rangle &=& \bigl\langle g_n\, ,\, 
ih_n^{-1}[P_\chi (h_n),a_{h_n}^w]g_n\bigr\rangle \, +\, O(h_n)\\
&=& \bigl\langle g_n\, ,\, (\{p_\chi ,a\})_{h_n}^wg_n\bigr\rangle \, +\, O(h_n)\period 
\label{mes-g_n} 
\end{eqnarray}
By (\ref{caract-measure}), the r.h.s. of (\ref{mes-g_n}) goes to $\mu (\{p_\chi ,a\})$, as 
$n\to \infty$. As in \cite{jec2}, we replace $P(h_n)$ by $P(h_n)-\lambda _n$ in the commutator on the 
l.h.s. of (\ref{comm-expectation}) and expand the commutator. Using (\ref{L2-bounds}), we show that 
$a_n=o(1)$, as $n\to \infty$, yielding $\mu (\{p,a\})=0$.\\
{\bf 2)}
The second assertion was established in the proof of
Lemma~\ref{g_n}.\\ 
{\bf 3)}
Let $\tau\in C^\infty(\R ^d)$ with support
in $\hat{M}$. Since $\support g_n\subset \{|x|\leq R_1'\}$,   
\[\bigl|\bigl\langle \tau ^2 g_n\, ,\, h_n^2\Delta _xg_n\bigr\rangle \bigr|
\ \leq \ |\langle \tau ^2g_n\, ,\, (P(h_n)-\lambda)g_n\rangle | \, +\, O(n^0)\comma \]
where $O(n^0)$ means $O(1)$ as $n\to\infty$. Thus
\begin{eqnarray*}
\bigl\langle ih_n\nabla_xg_n\, ,\, \tau ^2 ih_n\nabla _xg_n\bigr\rangle 
&\leq& 2h_n\, \bigl|\bigl\langle (\nabla _x\tau )g_n\, ,\, \tau ih_n\nabla _xg_n\bigr\rangle \bigr|
\, +\, O(n^0)\\
\|\tau ih_n\nabla _xg_n\|^2&\leq& O(h_n)\cdot \|\tau ih_n\nabla _xg_n\|\, +\, O(n^0)\comma 
\end{eqnarray*}
yielding the boundedness of $(\|\tau ih_n\nabla g_n\|)_n$. \qed

We introduce 
\begin{equation}\label{traj-bornee}
B_{\pm}(\lambda)\ := \ \bigl\{\ph \in p^{-1}(\lambda);\ 0\leq \pm t\, \donne\,  
\pi _x\phi (t;\ph )\hspace{.3cm} \mbox{is bounded}\bigr\}
\end{equation}
and $B(\lambda):=B_{+}(\lambda)\cap B_{-}(\lambda)$. By
(\ref{result-kn}), 
the non-trapping condition (\ref{non-trapping}) exactly means that $B_{+}(\lambda)$ and 
$B_{-}(\lambda)$ are empty. 
\begin{proposition}\label{supp-maj}
Let $d\geq 3$ if $N=0$ else let $d=3$. The measure $\mu$ is nonzero.\\
If $N=0$, $\mu$ vanishes near the (repulsive) singularities, 
is invariant under the complete flow $t\donne \phi ^t$, and 
$\support \mu \subset B(\lambda)$.\\ 
If $N>0$, then, outside the attractive 
singularities, $\mu$ is supported in $B(\lambda)$ that is 
\begin{equation}\label{supp-mu-B}
\support \mu \cap T^\ast (\R^3\setminus \cS)\ \subset \ B(\lambda)\period
\end{equation}
\end{proposition}
\Pf For the case of purely repulsive singularities (i.e.\ $N=0$)
the proof is given in 
Subsection~\ref{repuls-sing}. The other case appears in 
Subsection~\ref{general-case}. \qed

\begin{remark}\label{resonant-state}
If (\ref{size-resolvent}) is really false, one expects that the $f_n$ are
``close to some resonant state''. Proposition~\ref{supp-maj} and Proposition~\ref{almost-invariant}
 below roughly say that this resonant state should be microlocalized on
trajectories in
$B(\lambda)$. However, it does not give any information above 
the attractive singularities. If the potential $V$ is smooth (i.e.\ $N=N'=0$), 
the arguments used in \cite{jec2} actually prove Proposition~\ref{supp-maj} in this case.  
\end{remark}
\begin{lemma}\label{non-trapping-implies}
Let $d\geq 3$ if $N=0$ else let $d=3$. If 
$p$ is non-trapping at energy $\lambda$ (cf. (\ref{non-trapping}))
then $\mu=0$. 
\end{lemma}
\Pf Let $N=0$. By Proposition~\ref{supp-maj}, $\support \mu \subset B(\lambda)$, 
which is empty by the non-trapping condition. Thus $\mu =0$. The other case is 
treated in Subsection~\ref{general-case}. \qed

Now Proposition~\ref{supp-maj} and Lemma~\ref{non-trapping-implies} 
produce the desired contradiction.

\subsection{Repulsive singularities.}
\label{repuls-sing}

We show Proposition~\ref{supp-maj} for the case $N=0$, $d\geq 3$, by
first showing a decay estimate for the Fourier transform of the $g_n$'s.

Since we only have repulsive singularities, there exists some 
positive $c$ such that  
\begin{equation}\label{general-bound}
\bigl\langle g_n\, ,\, (-h_n^2\Delta _x)\,  g_n\bigr\rangle \, +\, \sum_{j=1}^{N'}
\bigl\langle g_n\, ,\, (1/|\cdot -s_j|)\,  g_n\bigr\rangle \ \leq \ c\, \langle g_n\, ,\, 
(P(h_n)-\lambda)g_n\rangle \, +\, O(n^0)\period 
\end{equation}
By Lemma~\ref{g_n}, $\|(P(h_n)-\lambda)g_n\|\to 0$ and the r.h.s of (\ref{general-bound}) 
is bounded. Now, we show that $\mu\neq 0$. Let $\chi\in C_0^\infty(\R ^d;\R)$ such that 
$0\leq \chi \leq 1$ and $\chi =1$ near $0$. Let us denote by $\cF g$ 
the Fourier transform of $g$. Setting $\chi_R(\xi)=\chi (\xi /R)$, for $R>0$ and $\xi\in \R^d$, we 
observe that 
\[\Bigl\langle \cF g_n\, ,\, \frac{(1-\chi_R)(h_n\cdot )}{|h_n\cdot |^2}\, 
|h_n\cdot |^2\cF g_n\Bigr\rangle \ \leq  \ \frac{O(n^0)}{R^2}\, \langle g_n\, ,\, (-h_n^2\Delta)
g_n\rangle \period \]
The bracket on the r.h.s is bounded uniformly w.r.t. $R$. Thus 
\begin{equation}\label{limlimsup}
\lim_R \limsup_n\langle \cF g_n\, ,\, (1-\chi_R)(h_n\cdot )\cF g_n\rangle 
\ = \ 0\period 
\end{equation}
Recall that, for all $n$, $\support g_n\subset \{|x|\leq R_1'\}$ (cf. Lemma~\ref{g_n}). 
By Proposition~\ref{total-mass}, this implies that $\|g_n\|^2\to \mu (\un )$, yielding 
$\mu \neq 0$. Now let $\tau\in C_0^\infty(\R;\R^+)$ be supported on a neighborhood of the
singularities such that $\tau =1$ near them. Since $V-\lambda$ is large and positive near 
the singularities, we can choose the support of $\tau$ such that, 
\begin{equation}\label{local-bound-potential}
\|\tau g_n\|^2\ \leq \ \langle \tau g_n\, ,\, (V-\lambda)\tau g_n\rangle \period
\end{equation}
\begin{equation}\label{local-bound}
\mbox{Thus}\hspace{.5cm}\bigl\langle \tau g_n\, ,\, (-h_n^2\Delta _x)\,  \tau g_n\bigr\rangle \, +\, \|\tau g_n\|^2
\ \leq \ 2\, \langle \tau g_n\, ,\, (P-\lambda)\tau g_n\rangle \, =\, o(1)\comma 
\end{equation}
using Lemma~\ref{basic-properties}. In particular, $\|\tau g_n\|^2\to \mu (\tau ^2)$ 
(cf. Proposition~\ref{total-mass}) and $\|\tau g_n\|\to 0$. 
Thus $\mu$ is supported away from the (repulsive) singularities. By Lemma~\ref{basic-properties}, 
we conclude that $\mu$ is invariant under the flow $(\phi ^t)_{t\in\R}$. If the trajectory 
$t\donne \pi _x\phi ^t(x,\xi)$ goes to infinity when $\pm t\to +\infty$, then the invariance 
of $\mu$ under the flow implies that $\mu$ vanishes on this trajectory. This shows that 
$\support \mu \subset B(\lambda)$ and finishes the proof of
Proposition~\ref{supp-maj} in the case $N=0$ and $d\geq 3$.

\subsection{The general case in dimension $3$.}
\label{general-case}

In this subsection, we assume that $N>0$ and $d=3$ and we give successively the 
proofs of Proposition~\ref{supp-maj} and Lemma~\ref{non-trapping-implies} 
(at the end of the subsection). In view of (\ref{limlimsup}) and of Proposition~\ref{total-mass}, 
we want to show that $\langle g_n,(-h_n^2\Delta _x)g_n\rangle$ is bounded 
to get $\mu\neq 0$. We also need a kind of invariance 
of $\mu$ under the pseudo-flow $\phi ^t$ (cf. (\ref{phi^t})). To realize this programme, we want to use the KS-transform (\ref{Lambda-ast}) to lift the property (\ref{L2-bounds}) in $\R^4$, locally near each attractive singularity. 

Let $(\tau_j)_{0\leq j\leq N}\in (C_0^\infty(\R^d;\R^+))^{N+1}$ be
such that \\
$\bullet$ $\sum _{j=0}^N\tau_j^2=1$ near $\{x\in \R^d; |x|\leq R_1'\}$,\\
$\bullet$ for $1\leq j\leq N$, $\tau_j=1$ near $s_j$ and is supported 
away from the other singularities, \\
$\bullet$ $\tau_0=1$ near the set of repulsive singularities and
is supported away from the other singularities. 

There exists $c>0$ such that 
\begin{eqnarray}\label{repuls-bound}
\sum_{j=N+1}^{N'}\bigl\langle \tau_0g_n\, ,\, (1/|x-s_j|)\,
\tau_0g_n\bigr\rangle \ \leq  c\, 
\langle \tau_0g_n\, ,\, 
(V-\lambda)\tau_0g_n\rangle \period \\
\mbox{Thus}\hspace{1.5cm}\bigl\langle \tau_0g_n\, ,\, (-h_n^2\Delta _x)\,  \tau_0g_n\bigr\rangle \, +\, 
\sum_{j=N+1}^{N'}\bigl\langle \tau_0g_n\, ,\, (1/|x-s_j|)\,  \tau_0g_n\bigr\rangle 
\nonumber\\
\leq  (1+c)\, \langle \tau _0g_n\, ,\, 
(P-\lambda)\tau _0g_n\rangle \, +\, O(n^0) \ = \ O(n^0)\period\label{repuls-energy-bound}
\end{eqnarray}
Here we used the fact that $\langle \tau_0g_n, (P-\lambda )
\tau_0g_n\rangle \to 0$, by Lemma~\ref{g_n} and Lemma~\ref{basic-properties}. 
Let $1\leq j\leq N$. For the same reason, $\langle \tau_jg_n, (P-\lambda )
\tau_jg_n\rangle \to 0$. Thus, since $(V-f_j/|\cdot -s_j|)\tau_j$ is bounded, 
\begin{equation}\label{leading-term}
\bigl|\bigl\langle \tau_jg_n\, ,\, (-h_n^2\Delta _x)\, \tau_jg_n\bigr\rangle \, +\, 
\bigl\langle \tau_jg_n\, ,\, (f_j/|\cdot -s_j|)\, \tau_jg_n\bigr\rangle\bigr| \ = \ O(n^0)\period
\end{equation}
We introduce the KS-transformation (cf. (\ref{Lambda-ast})) which is adapted to the singularity at $s_j$: 
$x=\cK _j (z_j):=s_j+\cK (z_j)$  (cf.\ (\ref{def-KS})) and, for $x\neq s_j$, 
\begin{equation}\label{Lambda-ast-s_j}
(x,\xi )\ =\ \cK _j^\ast (z;\zeta )\ :=\ (s_j,0)\, +\, \cK ^\ast (z;\zeta )\period 
\end{equation}
For all $n$, let $\tilde{g}_{n,j}:= g_n\rond \cK_j$. 
Let $\chi_j\in C_0^\infty(\R^3)$ such that $\chi_j\tau_j=\chi_j$, 
$\chi_j=1$ near $s_j,$ and $\chi_j\tau_k=0$, for $k\neq j$. Denote by
$\tilde{\chi}_j$ the function $\chi_j\rond \cK_j$. For $\lambda'\in \R$, we introduce 
the differential operator in $\R _{z_j}^4$ 
\begin{equation}\label{tilde-Pj}
\tilde{P}_j(h;\lambda ')\ := \ -h^2\Delta _{z_j}\, +\, \bigl((\tau_jV)\rond \cK_j-\lambda '\bigr)\, 
|\cdot |^2\comma
\end{equation}
which can be seen as the Weyl $h$-quantization of the symbol
\begin{equation}\label{tilde-pj}
T ^\ast\R^4\ni (z_j,\zeta_j)\ \donne \ \tilde{p}_{j,\lambda '}(z_j,\zeta_j)\ := \ 
|\zeta_j |^2\, +\, 
\bigl((\tau_jV)(\cK (z_j))-\lambda '\bigr)\cdot |z_j|^2\period
\end{equation}
Notice that $\tilde{p}_{j,\lambda '}\in \Sigma _{2;2}$. We can write, for 
$z_j\in \support\tilde{\chi}_j$, 
\begin{eqnarray}
|z_j|^2\bigl((\tau_jV)(\cK_j (z_j))-\lambda '\bigr) \!\!\! &=&\!\!\! f_j(s_j)\, + \, 
\bigl(f_j(\cK_j (z_j))-f_j(\cK_j (0))\bigr)+ \, |z_j|^2\bigl(W_j(\cK_j (z_j))-\lambda '\bigr)\nonumber\\
&=:& f_j(s_j)\, + \, \tilde{W}_{j,\lambda '}(z_j)\period \label{expression-potential}
\end{eqnarray}
So $\tilde{W}_{j,\lambda '}$ is a quadratic perturbation of the constant
$f_j(s_j)$, vanishing at $s_j$, and
\beq
\tilde{P}_j(h;\lambda ')=  -h^2\Delta _{z_j} +f_j(s_j)\, + \,
\tilde{W}_{j,\lambda '}.
\Leq{tilde-Pj:2}
\begin{lemma}\label{tilde-g_n}
Let $1\leq j\leq N$. The sequence $(\tilde{g}_{n,j})_n=(g_n\rond
\cK_j)_n$ is bounded in $\rL^2(\R^4)$. Up to subsequence, we may 
assume that it has a unique semiclassical measure $\tilde{\mu}_j$. 
Besides, 
\begin{eqnarray}
\forall n\in \N\comma &&\support \tilde{g}_{n,j}\ \subset \ \{z_j; |z_j|\leq (R_0+R_1')^{1/2}\}
\comma \label{support-tilde-g_n} \\
\mbox{and}\hspace{.2cm}&&\tilde{\chi}_j\tilde{P}_j(h_n;\lambda _n)\tilde{g}_{n,j}\ =\
o(h_n)\mbox{ in }\rL^2(\R^4)\period\label{tilde-Pj-small}
\end{eqnarray}
Let $\tilde{\phi}_j^s:=\tilde{\phi}_j(s;\cdot )$ be the Hamiltonian flow associated to
$(t,\lambda ',z,\zeta )\donne \tilde{p}_{j,\lambda '}(z,\zeta )$ by \eqref{cauchy-pb-extended}.
Let $\tilde{b}\in C_0^\infty(T^\ast\R^4)$ and $T_b:=\{s>0; \forall t\in
[0;s]\comma (\tilde{b}\circ \tilde{\phi}_j^t)(1-\tilde{\chi}_j)=0\}$. Then, for $s\in T_b$, 
\begin{equation}
\tilde{\mu}_j(\tilde{b})\ =\ \tilde{\mu}_j\bigl(\tilde{b}\circ
\tilde{\phi}_j^s\bigr)
\period \label{tilde-local-invariance}
\end{equation} 
\end{lemma}
\Pf 
$\bullet$
Eq. \eqref{support-tilde-g_n} follows from the scaling
$|\cK(z)|=|z|^2$
of the Hopf map (see \eqref{def-KS}) and
the estimate \eqref{L2-bounds} for the support of $g_n$ .\\ 
$\bullet$
Since $(-h_n^2\Delta_x+V-\lambda_n)g_n\ = \ o(h_n)$
and $g_n=O(n^0)$ in $\rL^2(\R^3)$, we use (\ref{def-KS}),
(\ref{change-KS}), and the arguments of Proposition 2.1 in \cite{gk}
to get 
\begin{equation}\label{tilde-L2-bounds}
|\cdot |^{-1}\tilde{\chi}_j\tilde{P}_j(h_n;\lambda
 _n)\tilde{g}_{n,j}\, =\,
 o(h_n)\hspace{.5cm}\mbox{and}\hspace{.5cm}|\cdot |\tilde{g}_{n,j}\,
 =\, O(n^0)\hspace{.5cm}\mbox{in}\hspace{.5cm}\rL^2(\R^4)\period 
\end{equation}
This yields (\ref{tilde-Pj-small}).\\ 
$\bullet$ 
Now, we show that
$\tilde{\chi}_j\tilde{g}_{n,j}=O(n^0)$ in $\rL^2(\R^4)$. Together with \eqref{tilde-L2-bounds}, this then will imply
the desired boundedness of $(\tilde{g}_{n,j})_n$ in $\rL^2(\R^4)$. \\
Thanks to (\ref{def-KS}), (\ref{change-KS}), and to Part 3 of Lemma~\ref{basic-properties}, 
\begin{eqnarray}
\|\un_{\support\nabla\tilde{\chi}_j}h_n\nabla _{z_j}\tilde{g}_{n,j}\| &=& 
O\bigl(\|\un_{\support\nabla\chi_j}h_n\nabla _xg_n\|\bigr)\ = \
O(n^0)\comma 
\label{esti-grad-away-from-sing}\\
\|\un_{\support\nabla\tilde{\chi}_j}\tilde{g}_{n,j}\|&=&
O\bigl(\|\un_{\support\nabla\chi_j}g_n\|\bigr)\ = \ O(n^0)\period\label{esti-away-from-sing}
\end{eqnarray}
Let $A_{n,j}\ :=\ (z_j\cdot h_n\nabla _{z_j} +h_n\nabla _{z_j}\cdot z_j)/(2i)$ and 
\[a_{n,j}\ :=\ \bigl\langle \tilde{g}_{n,j}\, ,\, h_n^{-1}\bigl[\tilde{P}_j(h_n;\lambda _n), 
i\tilde{\chi}_jA_{n,j}\tilde{\chi}_j\bigr]\, \tilde{g}_{n,j}\bigr\rangle\period \]
Expanding the commutator and using \eqref{tilde-L2-bounds}, we see, on one hand, that 
\begin{eqnarray}\label{comm-expan}
|a_{n,j}|&\leq & o(1)\cdot \bigl(h_n\|\un _{\support \nabla\tilde{\chi}_j}\tilde{g}_{n,j}\|\, +\, 
\|\tilde{\chi}_jih_n\nabla _{z_j}\tilde{g}_{n,j}\|\, +\, O(n^0)\bigr)\nonumber \\
&\leq & o(1)\cdot \|\tilde{\chi}_jih_n\nabla _{z_j}\tilde{g}_{n,j}\| \, +\, o(1)\comma 
\end{eqnarray}
thanks to (\ref{esti-away-from-sing}). On the other hand, writing
$2iA_{n,j}=2z_j\cdot h_n\nabla _{z_j}+4h_n$, 
\begin{eqnarray*}
a_{n,j}&=& \bigl\langle \tilde{g}_{n,j}\, ,\, 2\bigl[-h_n^2\Delta_{z_j},\tilde{\chi}_j^2\bigr]
\tilde{g}_{n,j}\bigr\rangle\, +\, \Im \bigl\langle z_j\cdot h_n\nabla _{z_j}\tilde{\chi}_j
\tilde{g}_{n,j}\, ,\, h_n^{-1}\bigl[-h_n^2\Delta_{z_j},\tilde{\chi}_j\bigr]\tilde{g}_{n,j}
\bigr\rangle\\
&& +\, \bigl\langle \tilde{\chi}_j\tilde{g}_{n,j}\, ,\, h_n^{-1}\bigl[\tilde{P}_j(h_n;\lambda _n), 
z_j\cdot h_n\nabla _{z_j}\bigr]\, 
\tilde{\chi}_j\tilde{g}_{n,j}\bigr\rangle\period 
\end{eqnarray*}
By \eqref{esti-grad-away-from-sing} and \eqref{esti-away-from-sing}, 
\begin{eqnarray*}
\bigl| a_{n,j}\, -\, \bigl\langle \tilde{\chi}_j\tilde{g}_{n,j}\, ,\, 
h_n^{-1}\bigl[\tilde{P}_j(h_n;\lambda _n), iA_{n,j}\bigr]\, 
\tilde{\chi}_j\tilde{g}_{n,j}\bigr\rangle \bigl| &\leq &
O(n^0)\cdot\|\tilde{\chi}_jih_n\nabla _{z_j}
\tilde{g}_{n,j}\| 
\, +\, O(n^0)
\period 
\end{eqnarray*}
As a differential operator, $h_n^{-1}[\tilde{P}_j(h_n;\lambda _n), iA_{n,j}]=
2(-h_n^2\Delta _{z_j})-z_j\cdot\nabla _{z_j}\tilde{W}_{j,\lambda _n}(z_j)$ 
(cf. (\ref{tilde-Pj:2})) and, by (\ref{support-tilde-g_n}), there exists some 
$c_j>0$ such that, for all $n$ and for all $z_j\in \support \tilde{g}_{n,j}$, 
$|z_j\cdot\nabla _{z_j}\tilde{W}_{j,\lambda _n}(z_j)|\leq
c_j|z_j|^2$. By \eqref{tilde-L2-bounds}, 
\begin{eqnarray*}
\bigl| a_{n,j}\, -\, \bigl\langle \tilde{\chi}_j\tilde{g}_{n,j}\, ,\, -2h_n^2\Delta_{z_j}\, 
\tilde{\chi}_j\tilde{g}_{n,j}\bigr\rangle\bigl| &\leq &
O(n^0)\cdot\|\tilde{\chi}_jih_n\nabla _{z_j}
\tilde{g}_{n,j}\| 
\, +\, O(n^0)\period
\end{eqnarray*}
This, together with (\ref{comm-expan}), implies that 
\begin{equation}\label{esti-laplacian}
0\ \leq \ \bigl\langle \tilde{\chi}_j\tilde{g}_{n,j}\, ,\, -2h_n^2\Delta_{z_j}\, 
\tilde{\chi}_j\tilde{g}_{n,j}\bigr\rangle\ \leq \
O(n^0)\cdot\|\tilde{\chi}_jih_n
\nabla _{z_j}\tilde{g}_{n,j}\| 
\, +\, O(n^0)\period
\end{equation}
Writing 
\begin{eqnarray*}
\bigl\langle \tilde{\chi}_j\tilde{g}_{n,j}\, ,\, -h_n^2\Delta_{z_j}\, 
\tilde{\chi}_j\tilde{g}_{n,j}\bigr\rangle &=&\|\tilde{\chi}_jih_n\nabla _{z_j}\tilde{g}_{n,j}\|^2\, +\,
h_n^2\|(\nabla _{z_j}\tilde{\chi}_j)\tilde{g}_{n,j}\|^2\\
&&\, +\, 2h_n\rRe \langle (\nabla _{z_j}\tilde{\chi}_j)\tilde{g}_{n,j}\, ,\, \tilde{\chi}_jih_n
\nabla _{z_j}\tilde{g}_{n,j}\rangle 
\end{eqnarray*}
and using again (\ref{esti-grad-away-from-sing}) and (\ref{esti-away-from-sing}), we arrive at 
\[\|\tilde{\chi}_jih_n\nabla _{z_j}\tilde{g}_{n,j}\|^2\ \leq \ O(n^0)\cdot\|\tilde{\chi}_jih_n
\nabla _{z_j}\tilde{g}_{n,j}\| 
\, +\, O(n^0)\period\]
This yields 
\begin{equation}\label{esti-grad-squared-laplacian}
\|\tilde{\chi}_jih_n\nabla _{z_j}\tilde{g}_{n,j}\|\ = \ O(n^0)\hspace{.5cm}
\mbox{and}\hspace{.5cm}\bigl\langle \tilde{\chi}_j\tilde{g}_{n,j}\, ,\, -h_n^2\Delta_{z_j}\, 
\tilde{\chi}_j\tilde{g}_{n,j}\bigr\rangle\ = \ O(n^0)\period
\end{equation}
Now
\begin{eqnarray*}
\bigl\langle \tilde{\chi}_j\tilde{g}_{n,j}\, ,\, \, \tilde{P}_j(h_n;\lambda _n)
\tilde{\chi}_j\tilde{g}_{n,j}\bigr\rangle &=&\bigl\langle \tilde{\chi}_j\tilde{g}_{n,j}\, ,\, \, 
\tilde{\chi}_j\tilde{P}_j(h_n;\lambda _n)\tilde{g}_{n,j}\bigr\rangle \, +\, \bigl\langle
\tilde{\chi}_j\tilde{g}_{n,j}\, ,\, \, \bigl[-h_n^2\Delta_{z_j},\tilde{\chi}_j\bigr]
\tilde{g}_{n,j}\bigr\rangle 
\end{eqnarray*}
and is bounded by \eqref{tilde-L2-bounds}, (\ref{esti-grad-away-from-sing}), and (\ref{esti-away-from-sing}). Thus 
\beq
\bigl\langle \tilde{\chi}_j\tilde{g}_{n,j}\, ,\,
-h_n^2\Delta_{z_j}\, \tilde{\chi}_j\tilde{g}_{n,j}
\bigr\rangle \, +\,
f(s_j)\|\tilde{\chi}_j\tilde{g}_{n,j}\|^2\, +\, \langle\tilde{\chi}_j\tilde{g}_{n,j}\, ,\, 
\tilde{W}_{j,\lambda _n}\tilde{\chi}_j\tilde{g}_{n,j}\rangle \ = \ O(n^0)\period
\Leq{OOOn0}
In (\ref{OOOn0}), the first and third terms are $O(n^0)$, by (\ref{esti-grad-squared-laplacian}) and by \eqref{tilde-L2-bounds} respectively. 
Since  $f_j(s_j)\neq 0$, we conclude that 
$(\tilde{\chi}_j\tilde{g}_{n,j})_n$ is bounded in $\rL^2(\R^4)$. \\
$\bullet$ We now show the invariance \eqref{tilde-local-invariance}. 
It suffices to show that, for all $\lambda\in \R$ and all $\tilde{b}\in
C_0^\infty(T^\ast\R^4)$ such that $\tilde{b}(1-\tilde{\chi}_j)=0$, 
$\tilde{\mu}_j(\{\tilde{p}_{j,\lambda}, \tilde{b}\})=0$. Take such a
$\tilde{b}$ and $\lambda\in \R$. Since 
$\tilde{b}_{h_n}^w$ is uniformly bounded, 
\[\bigl\langle \tilde{g}_{n,j} ,ih_n^{-1}\bigl[\tilde{\chi}_j\tilde{P}_j(h_n;\lambda _n), 
\tilde{b}_{h_n}^w\bigr]\tilde{g}_{n,j}\bigr\rangle \ = \ o(n^0)\comma \]
by expanding the commutator, using (\ref{tilde-Pj-small}), and using
the boundedness in $\rL^2(\R^4)$ of $(\tilde{g}_{n,j})_n$. 
Now we compute the leading term of the commutator and arrive at 
\begin{eqnarray*}
o(n^0)&=& \bigl\langle \tilde{g}_{n,j} ,\{\tilde{\chi}_j\tilde{p}_{j,\lambda_n}, \tilde{b}\}_{h_n}^w
\tilde{g}_{n,j}\bigr\rangle \, +\, O(h_n) 
\ =\ \bigl\langle \tilde{g}_{n,j} ,\{\tilde{\chi}_j\tilde{p}_{j,\lambda}, \tilde{b}\}_{h_n}^w
\tilde{g}_{n,j}\bigr\rangle \, +\, o(n^0)\\
&=&\bigl\langle \tilde{g}_{n,j} ,\{\tilde{p}_{j,\lambda}, \tilde{b}\}_{h_n}^w
\tilde{g}_{n,j}\bigr\rangle \, +\, o(n^0)\comma
\end{eqnarray*}
since $\tilde{\chi}_j=1$ on the support of $\tilde{b}$. Thus $\tilde{\mu}_j(\{\tilde{p}_{j,\lambda}, 
\tilde{b}\})=0$. \qed

\Pfof{Proposition~\ref{supp-maj}} Let $1\leq j\leq N$. The boundedness
of the sequence $(\tilde{\chi}_j\tilde{g}_{n,j})_n$ in $\rL^2(\R^4)$ precisely means that
$(\langle\chi_jg_n\, ,\, (1/|\cdot -s_j|)\chi_jg_n\rangle)_n$ is
bounded (cf. (\ref{change-KS})) 
and so is 
also $(\langle\tau_jg_n\, ,\, (1/|\cdot -s_j|)\tau_jg_n\rangle)_n$. By (\ref{leading-term}), 
this implies that $(\langle\tau_jg_n\, ,\, -h_n^2\Delta_x\tau_jg_n\rangle)_n$ is bounded. 
By the IMS localization formula (cf. Chapter 3.1 of \cite{cfks}), 
\begin{equation}\label{esti-grad-squared-laplacian-bis}
\langle g_n\, ,\, -h_n^2\Delta_x g_n\rangle \ = \ \sum_{j=0}^N\langle\tau_jg_n\, ,\, -h_n^2\Delta_x
\tau_jg_n\rangle \ -\
h_n^2\sum_{j=0}^N\|(\nabla _x \tau_j)g_n\|^2\ = \ O(n^0)\comma
\end{equation}
thanks to (\ref{repuls-energy-bound}). As in Subsection~\ref{repuls-sing}, we can 
derive (\ref{limlimsup}) and prove that $\mu\neq 0$. \\
Consider a trajectory $(\phi (t;\ph _0))_{t\not\in \J(\ph _0))}$ such
that $\pi_x\phi (t;\ph _0)$ goes to infinity as $t\to \pm \infty$. 
If it does hit a singularity then $\pi_x\phi (t;\ph _0)$ must come from
infinity, hit the singularity and then go back to infinity ($\J(\ph
_0)$ contains one point). Since $\mu$ vanishes on some $\{\ph \in T^\ast\R^3;|x|\geq C\}$, 
$\mu$ vanishes near the tail(s) of $(\phi (t;\ph _0))_{t\not\in \J(\ph _0))}$ which is
(are) inside this set. By invariance (cf. Lemma~\ref{basic-properties}), $\mu$ vanishes near each 
$\phi (t;\ph _0)$, for $t\not\in \J(\ph _0)$. This proves (\ref{supp-mu-B}). \qed

\Pfof{Lemma~\ref{non-trapping-implies}} 
Let $1\leq j\leq N$ and $\tilde{\tau}\in C_0^\infty(\R^4)$ with
$\tilde{\tau}(1-\tilde{\chi}_j)=0$ and $\tilde{\tau}=0$ near
$z_j=0$. Then $|\tilde{\tau}|^2\tilde{\mu}_j$
is the semiclassical measure of $(\tilde{\tau}\tilde{g}_{n,j})_n$ (see
\cite{gl}). We may assume that $\tau =\tilde{\tau}\circ J_{j,+}$ is well
defined. By \eqref{change-KS}, $\|\tau g_n|\cdot |^{-1/2}\|^2=\|\tilde{\tau}\tilde{g}_{n,j}\|^2$.
By \eqref{def-KS}, $\tau _1:=\tau |\cdot |^{-1/2}$ is smooth. Thus $\langle \tau_1g_n, (P-\lambda )
\tau_1g_n\rangle \to 0$, by Lemma~\ref{g_n} and
Lemma~\ref{basic-properties}. This yields the bound \eqref{esti-grad-squared-laplacian-bis}
and Eq. (\ref{limlimsup}) with $g_n$ replaced by $\tau _1g_n$. By
Proposition~\ref{total-mass}, $\|\tau _1g_n\|^2\to |\tau _1|^2\mu (\un
)$. But the latter is zero since, by Proposition~\ref{supp-maj} and the non-trapping 
assumption, $\mu$ may only have mass above the attractive
singularities. Thus $\lim \|\tilde{\tau}\tilde{g}_{n,j}\|=0$. This
implies that $\tilde{\chi}_j\tilde{\mu}_j$ may only have mass above
$z_j=0$. \\
Now $\tau \in C_0^\infty(\R^3)$ 
supported near $s_j$ and inside the set $\chi_j^{-1}(1)$. Let $\tilde{\tau}=
\tau\circ \cK _j$. Thanks to \eqref{def-KS}, we can choose the support of
$\tau$ small enough such that, for some $s'>0$, 
$\support \tilde{\tau}\circ \tilde{\phi}(s';\cdot )\subset\tilde{\chi}_j^{-1}(1)$ 
and $\tilde{\tau}\circ \tilde{\phi}(s';\cdot )=0$ near $\{0\}\times \R^4\subset T^\ast\R^4$.
By (\ref{change-KS}), $\|\tau g_n\|^2=\|\tilde{\tau}\tilde{g}_{n,j}|\cdot |\, \|^2$. 
By (\ref{esti-grad-squared-laplacian}), Eq. (\ref{limlimsup}) holds true with $g_n$ replaced 
by $\tilde{g}_{n,j}$. 
This implies, by Proposition~\ref{total-mass}, that $\|\tilde{\tau}\tilde{g}_{n,j}|\cdot |\, \|^2\to
\tilde{\mu}_j(|\cdot |^2\tilde{\tau}^2)$. By \eqref{tilde-local-invariance}, 
$\tilde{\mu}_j(|\cdot |^2\tilde{\tau}^2)=\tilde{\mu}_j((|\cdot
|^2\tilde{\tau}^2)\rond\tilde{\phi}(s';\cdot))=0$. 
Therefore $\lim \|\tau g_n\|^2=0$, yielding $\mu=0$ near $s_j$. Thus
$\mu =0$. \qed

Actually, if trapping occurs, we have the following stronger result
on the measure $\mu$. 

\begin{proposition}\label{almost-invariant}
Let $N>0$ and $d=3$. If $\ph \in \support \mu\cap T^\ast
(\R^3\setminus \cS)$ and $t\not\in \J(\ph )$ then $\phi (t;\ph)\in
\support \mu$. 
\end{proposition}

\Pf Let $1\leq j\leq N$. Let $\ph _0:=(x_0,\xi_0)\in p^{-1}(\lambda)$ such that 
$\chi_j=1$ near $x_0$. By the properties of the KS-transform (\ref{Lambda-ast-s_j}) 
(cf. (\ref{sol-zeta})), 
there exists $\phz _0=(z_0,\zeta_0)\in T^\ast\R ^4$ such that $\ph _0=\cK _j^\ast 
(\phz _0)$. Let $\pht _0=(0,p(\ph _0))=(0,\lambda)$. We consider the trajectory 
$\{\pi_x\phi ^t(\ph _0), t\in \R\}$ and assume that it hits the singularity $s_j$ at 
time $t_0$. Let $t'>t_0$ such that $\chi_j(\pi _x\phi (t';\ph _0))=1$. There exists some $s'\in\R$ such that $t'=t_j(s';\pht _0,\phz _0)$ (cf. (\ref{phi-phitilde})). Here $t_j(s;\pht ,\phz )$ is the first component of the flow $\tilde{\phi}_j(s;\pht ,\phz )$ given by \eqref{phi-phitilde} with $\tilde{p}$ replaced by \eqref{tilde-pj}. 
Let $\tau _0\in C_0^\infty(\R^3)$ such that  $\chi_j=1$ near $\support \tau_0$, $\tau _0=1$ near $x_0$, and $\tau _0=0$ near $s_j$. The semiclassical measure $\mu _1$ of the sequence $(\tau _0g_n)_n$, viewed as a bounded sequence in $\rL^2(\R^3\times S^1)$, is 
$\mu \otimes 1\otimes \delta _0$ on $T^\ast\R^3\times T^\ast S^1$. Let 
$\psi\in C_0^\infty(\R)$ such that $\psi =1$ near $0$ and $K_0\subset\subset\R^3$ 
be a vicinity of $\xi _0$. Let $a\in 
C_0^\infty(T^\ast\R^3)$ such that $\tau _0=1$ near $\pi _x\support a$ and $\pi _\xi 
\support a\subset K_0$. For $(\ph ;\theta
^\ast):=(x;\xi ;\theta ;\sigma)\in T^\ast\R^3\times T^\ast S^1$, set $a_1(\ph ;\theta
^\ast)=\psi (\sigma )a(\ph)$. Let $\psi _1+\psi _2=1$ be a smooth
partition of unity on $S^1$. Notice that 
\begin{equation}\label{mu_1}
\mu (a)\ =\ \tau _0\mu (a)\ =\ \tau _0\mu _1(a_1)\ =\ \sum
_{k=1}^2\tau _0\mu _1(a_1\psi _k)\period 
\end{equation}
For each $k\in \{1;2\}$, we may apply Proposition~\ref{diffeo} with 
$u_n=\tau _0g_n\psi_k\in\rL^2(\R^3\times S^1)$ and $\Phi =(\cK_j ,\cA _{j,+})$, since 
$(\cK_j ,\cA _{j,+})$ is a local diffeomorphism near $\support \tau_0\times
\support\psi _k$ by \eqref{local-diffeo}. Thus $((\tau_0\psi _k)\circ (\cK_j ,\cA _{j,+}))
\tilde{\mu}_j(\tilde{b}_k)=\tau _0\mu _1(a_1\psi _k)$, where $\tilde{b}_k=(a_1\psi _k)\circ (\cK_j ,\cA _{j,+})_c$, since $((\tau_0\psi _k)\circ (\cK_j ,\cA _{j,+}))\tilde{\mu}_j$ is the semiclassical measure
of $((\tau _0g_n\psi _k)\circ (\cK_j ,\cA _{j,+}))_n$. Now we can choose $K_0$ and $\support\tau_0$ small enough such 
that, for all $k\in \{1;2\}$, $\tilde{b}_k\circ \tilde{\phi}_j^{s'}=0$ near 
$\{0\}\times \R^4$ and $(1-\tilde{\chi}_j)\tilde{b}_k\circ
\tilde{\phi}_j^{t}=0$, for $0\leq t\leq s'$. Thus
\eqref{tilde-local-invariance} holds true with $s=s'$ and
$\tilde{b}=\tilde{b}_k$. Let $\tilde{\tau}_k\in C_0^\infty(\R^4)$ such
that  $\tilde{\chi}_j=1$ near $\support \tilde{\tau}_k$,
$\tilde{\tau}_k=1$ near $\pi _z\support \tilde{b}_k\circ \tilde{\phi}_j^{s'}$, and $\tilde{\tau}_k=0$ 
near $z_j=0$. We may assume that $\cJ _{j,+}$ is a local
diffeomorphism with local inverse $(\cK_j ,\cA _{j,+})$
near $\pi _z\support \tilde{b}_k\circ \tilde{\phi}_j^{s'}$
(cf. \eqref{local-diffeo}). Thus we can apply Proposition~\ref{diffeo} with 
$u_n=\tilde{\tau}_k\tilde{g}_{n,j}\in\rL^2(\R^4)$ and $\Phi =\cJ
_{j,+}$. This yields $\tilde{\tau}_k\tilde{\mu}_j(\tilde{b}_k\circ
\tilde{\phi}_j^{s'})=\tau _k\mu _1(a_{s',k})$, where $\tau
_k=\tilde{\tau}_k\circ \cJ _{j,+}$ and $a_{s',k}=\tilde{b}_k\circ
\tilde{\phi}_j^{s'}\circ (\cJ _{j,+})_c$. Now we see that, if 
$\mu$ is zero near $\phi (t';\ph _0)$, then we can choose 
$K_0$ and $\support\tau_0$ small enough such that $\tau _k\mu
_1(a_{s',k})=0$, for $k\in \{1;2\}$. By \eqref{mu_1}, this implies that $\mu (a)=0$,
for $a$ with small enough support near $\ph _0$. Since we can 
reverse the time direction, we get the desired result.  \qed

\subsection{A simpler proof for weighted $\rL^2$ estimates.}
\label{simpler-case}

In Subsections~\ref{main-lines}, \ref{repuls-sing}, and~\ref{general-case}, we proved that the 
non-trapping condition implies the Besov estimate (\ref{size-resolvent}). By (\ref{order}), the latter 
implies the existence of some $C>0$ such that, for all $s>1/2$, 
\begin{equation}\label{bound-weighted-L2}
\sup_{\stackrel{\Re z\in I}{\Im z\neq 0}}\, \|R(z;h)\|_{\rL _s^2,\rL _{-s}^2}\ \leq \ C\cdot 
h^{-1}\comma
\end{equation}
a weighted $\rL^2$ estimate. This derivation of (\ref{bound-weighted-L2}) from the non-trapping 
condition uses Proposition~\ref{compact-local}, the proof of which is based on arguments borrowed 
from \cite{cj}. 
The latter are rather involved since, in \cite{cj}, the potential is assumed to be $C^2$ only. 
In particular, a special pseudodifferential calculus, adapted to this low regularity, is used there. 
Since our potential here is $C^\infty$ outside the singularities, we want to give 
a simpler proof of the following, slightly weaker result.
\begin{proposition}\label{easier-result}
Under the assumptions of Theorem~\ref{main}, we assume that $p$ is non-trapping at each energy 
$\lambda\in I_0$. Then, for any compact interval $I\subset I_0$ and any $s>1/2$, there exists
$C_s>0$ such that (\ref{bound-weighted-L2}) holds true with $C=C_s$.
\end{proposition}
\Pf Let $d\geq 3$. We can follow the arguments in Subsections~\ref{main-lines}, \ref{repuls-sing}, 
and~\ref{general-case}, if Proposition~\ref{compact-local} is replaced by 
\begin{equation}\label{infinity-weighted-L2}
\exists R_0'>R_0\, ;\hspace{.4cm} \lim_{n\to\infty}\bigl\|\un_{\{|\cdot |>R_0'\}}f_n
\bigr\|_{\rL^2_{-s}}\ =\ 0\period 
\end{equation}
Indeed, for functions localized in $\{x\in \R^d; |x|\leq R_0'\}$, the norms $\|\cdot \|_B$ and 
$\|\cdot \|_{\rL _s^2}$ are equivalent and so are the norms $\|\cdot\|_{B^\ast}$ and 
$\|\cdot \|_{\rL _{-s}^2}$. 
So we are left with the proof of (\ref{infinity-weighted-L2}). We follow 
the proof of Proposition~\ref{compact-local} and arrive at (\ref{regular-terms}). Now, 
by \cite{jec2}, we can find $c>0$, a function $\chi_1\in C_0^\infty(\R ^d)$, and a 
symbol $a\in\Sigma _{0;0}$ satisfying the following properties. 
The function $\chi _1=1$ on a large enough neighbourhood of $0$ and of 
the support of $\chi $ and 
\[\alpha _n \ \geq \ c\cdot \bigl\|(1-\chi _1)\theta (P_\chi (h_n))f_n
\bigr\|_{\rL^2_{-s}}^2\, +\, o(1)\period \]
Following again the proof of Proposition~\ref{compact-local}, we get (\ref{infinity-weighted-L2}). \qed

\section{On the validity of the non-trapping condition.}
\label{when-non-trapping?}
\setcounter{equation}{0}

The aim of this last section is to provide examples both of validity and of invalidity of the non-trapping 
condition (\ref{non-trapping}). As we shall see in Corollary~\ref{effective-trapping} below, the 
non-trapping property is seldom fulfilled if there is some singularity ($N'>0$),
even at positive energies. This is in strong contrast to the smooth case, for which $p$ is always 
non-trapping at large enough positive energies. 

To study the non-trapping condition (\ref{non-trapping}) when an attractive singularity is present 
(and $d=3$), we need to review the regularization of the Hamilton flow of $p$, described 
in Section~\ref{pseudo-sing}, in a more sophisticated way. Recall that 
$\hat{M}=\R^3\setminus\cS$. Let $\omega_0$ be the natural symplectic two-form on 
$T^\ast\R ^3$ given by $\sum_{j=1}^3dx_i\wedge d\xi_i$ and also its restriction to $\hat{P}$. 
It is well known (see \cite{kn}, Thm.\ 5.1) that there exists an extension $(M,\omega ,m)$ 
of the Hamiltonian system $(\hat{P}, \omega_0, p)$, where as a set 
the six-dimensional smooth manifold $M$ equals
\[M\ := \ \hat{P}\, \cup\, \bigcup_{i=1}^N
(\R\times S^2)\period \]
Here the $i$th copy $\R\times S^2$ parameterizes energy and direction of the particle 
colliding with the attractive singularity $s_i$. Using the symplectic form $\omega$ on $M$, 
the Hamiltonian function $m\in C^\infty(M)$ generates a smooth complete flow 
\begin{equation}\label{flow-above}
\Phi : \R\times M\dans M \hspace{.4cm}\comma\hspace{.4cm}(t;\ph )
\donne \Phi (t;\ph )=: \Phi ^t(\ph )\period
\end{equation}
A collision time for $\ph \in M$ is a time $t_0$ such that $\Phi(t_0; \ph )\not\in 
\hat{P}$. If $t$ is not a collision time for 
$\ph\in \hat{P}$ then $\Phi(t; \ph )=\phi (t;\ph )$, 
defined just before (\ref{collision-time}). 
\begin{proposition}\label{trapping}
Consider for $d=2$ or $3$ a regular value $\lambda>0$ of $V$. 
If the set
\[\cH_\lambda \ :=\ \bigl\{x\in\R ^d ;\,  V(x)\geq\lambda\ 
\mbox{or}\ x\in\cS \bigr\}\]
is not homeomorphic to a $d$-dimensional ball or a point, then $p$ is
trapping at energy $\lambda$, i.e.\ (\ref{non-trapping}) is false.
\end{proposition}
\Pf We write $\cH_\lambda$ as the disjoint union $\tilde{\cH}_\lambda\dot{\cup}\{s_1,\ldots,
s_N\}$ with
\[\tilde{\cH}_\lambda \ :=\ \bigl\{x\in\R ^d ;\,  V(x)\geq\lambda\ \mbox{or}\ x\in\{s_{N+1},
\ldots , s_{N'}\}\bigr\}\period \]
Then $\tilde{\cH}_\lambda$ is a $d$--dimensional manifold with boundary, since by 
assumption $\lambda$ is a regular value of $V$. It is compact since 
by assumption $\lim_{|x|\to\infty}V(x)=0$ but $\lambda>0$. Notice that $\tilde{\cH}_\lambda$ 
is a neighbourhood of the repulsive singularities $s_{N+1},\ldots,s_{N'}$, but there exist 
neighbourhoods of the attractive singularities $s_1,\ldots,s_N$ that are disjoint from $\tilde{\cH}_\lambda$.
In the presence of repulsive singularities $\tilde{\cH}_\lambda$ is nonempty. In any
case, $\cH _\lambda$ is a nonempty compact set. We denote by ${\rm Int}({\cH}_\lambda)$ the interior of 
$\cH _\lambda$. Now we assume that $\cH_\lambda$ is not homeomorphic 
to a $d$-dimensional ball nor to a point, and we construct 
a periodic orbit, thus proving trapping. We discern two cases. \\
\underline{First case}: $\cH_\lambda$ has two or more connected components.\\ 
Here the idea is to construct a periodic orbit (using curve shortening), whose projection on 
configuration space is a curve connecting two components of $\cH_\lambda$. Let
$g_{\rm Euclid}$ denotes the euclidean metric on $\R^d$. 
We now use the Jacobi metric $\hat{g}_\lambda$, defined on
$\R^d\setminus\cH_\lambda$ by 
\begin{equation}\label{jac}
\hat{g}_\lambda(q):=(\lambda-V(q))g_{\rm Euclid}.
\end{equation}
It is known (see e.g. \cite{kk} and \cite{bn}) that for regular curves
$c:[0,1]\to\R^d\setminus{\rm Int}({\cH}_\lambda)$ with $c(1)=s_i\ (i\leq N)$ 
the length 
\[\cL(c):=\lim_{t\nearrow1}
\int_0^t\sqrt{\hat{g}_\lambda(c(s))\big(\dot{c}(s),\dot{c}(s)\big)}ds\] 
is finite. By compactness of 
$\cH_\lambda$ the number $\ell$ of connected components
of $\cH_\lambda$ is finite. Denoting them by 
$\cH_{\lambda ;1},\ldots,\cH_{\lambda ;\ell}$, for $1\leq i<j\leq \ell$ 
\[D_\lambda(i,j):=\inf_{c:c(0)\in\cH_{\lambda ;i}\, ,\, c(1)\in\cH_{\lambda ;j}}\cL(c)>0,\]
that is, the different components have positive geodesic distances.\\
Taking $R$ large enough, we can ensure that these mutual distances
are smaller than the corresponding geodesic distance of the $\cH_{\lambda ;i}$ to 
the region $\{x\in\R^d; |x|\geq R\}$.\\
The (standard) approach is to consider the negative gradient flow of the {\em energy 
functional}
\[\cE(c):=\int_0^1\hat{g}_\lambda(c(s))\big(\dot{c}(s),\dot{c}(s)\big)\,ds,\qmbox{with} 
c(0)\in\cH_{\lambda ;i_0}\mbox{ and }
c(1)\in\cH_{\lambda ;i_1}\]
in order to approximate geodesic segments, which are then critical points of $\cE$ with respect to 
these boundary conditions.\\
Due  to the degeneracy of the Jacobi metric (\ref{jac}) at $\partial (\R^d\setminus\cH_\lambda)$
still no Palais--Smale condition is satisfied for $\cE$, 
that is, a vanishing gradient of $\cE$ at $c$
does not ensure that $c$ is a geodesic (see Klingenberg \cite{kl}, 
Chapter 2.4 for a discussion of the Palais--Smale condition).\\
However, as $\lambda$ is assumed to be a regular value of $V$, the regularization 
technique devised by Seifert in \cite{s} and later by Gluck and Ziller in 
\cite{gz} can be applied to yield a geodesic segment of length equal to 
$D_\lambda(i_0,i_1)=\min_{i<j}D_\lambda(i,j)>0$, with $c(0)\in\cH_{\lambda ;i_0}$ and 
$c(1)\in\cH_{\lambda ;i_1}$.\\
We denote the restriction of the 
flow $\Phi^t$ to $m^{-1}(\lambda)$ by $\Phi_\lambda^t$. 
Away from the end points, and 
up to time parameterization, the 
geodesic segment in the Jacobi metric corresponds to a segment of a $\Phi_\lambda^t$--solution curve.
See \cite{am}, Thm. 3.7.7 for a proof.\\
This segment is part of a periodic orbit, whose period is twice
the time needed to parametrize the segment:
\begin{itemize}
\item
If $V(c(i_k))=\lambda$, then (by our regularity assumption for the value $\lambda$)
$\nabla V(c(i_k))\neq0$.
Furthermore the geodesic segment at this point has a normalized tangent 
\[\lim_{s\nearrow i_k} \|c(s)-c(i_k)\|^{-1}(c(s)-c(i_k))\]
which is parallel to $\nabla V(c(i_k))$ (see \cite{gz}, Sect. 6). 
Thus the solution curve can be continued by {\em time reversal} (cf. (\ref{time-reversal}))
\begin{equation}\label{tire}
c(i_k+s):=c(i_k-s)\qquad(s>0). 
\end{equation}
\item
Similarly, if instead 
$c(i_k)\in\{s_1,\ldots,s_N\}$, that is, $c(t)$ converges to an attracting singularity,
then, time reversal (\ref{tire}) again continues the 
geodesic segment $c$ and thus the 
$\Phi_\lambda^t$--solution curve as well.
\end{itemize}
In both cases we thus constructed a periodic $\Phi_\lambda$-orbit.\\
\underline{Second case}: $\cH_\lambda$ has only one component, which however is not 
homeomorphic to a $d$-dimensional ball nor a point.
Thus it is a connected compact $d$--dimensional submanifold of $\R^d$ with
boundary not homeomorphic to ${\mathbb S}^{d-1}$. 
\begin{itemize}
\item
If $\R^d\setminus{\rm Int}({\cH}_\lambda)$ contains a compact connected component, then this arises 
as the projection on configuration space of a connected component
of the regularized energy surface $m^{-1}(\lambda)$.
This flow-invariant component is compact too, and 
thus consists of trapped orbits.
\item
If, however $\R^d\setminus{\rm Int}({\cH}_\lambda)$ does not contain a compact
component, it necessarily is connected since $d\geq2$ and $\cH_\lambda$ is
compact. In this situation, the boundary $\partial\cH_\lambda$ consists of one component, 
which is not homeomorphic to ${\mathbb S}^{d-1}$. In this situation Corollary 3.3 of \cite{kn1} ensures 
the existence of a periodic so--called brake orbit, that is a 
trapped orbit in the terminology of our paper (although \cite{kn1} treats smooth
potentials, in the case at hand all singularities of our potential are repelling.
Thus the dynamics at energy $\lambda$ is unaffected by the singularities.). 
\qed
\end{itemize}

A converse of Proposition~\ref{trapping} does not hold true in general. That is, there are 
potentials like Yukawa's potential $V(x)=-e^{-|x|}/|x|$ for which $\cH_\lambda$ consists only 
of one point but still there are trapped orbits for small $\lambda>0$, see 
\cite{kk}. Yet Proposition~\ref{trapping} gives us the 
\begin{corollary}\label{effective-trapping}
Consider for $d=2$ or $3$ a regular value $\lambda>0$ of $\,V$.\\ 
If $N>1$ or if $N=1$ and $N'>N$, then $p$ is trapping at 
energy $\lambda$.\\
If $N'\geq 2$ then $p$ is trapping at 
energy $\lambda$, for $\lambda$ large enough. 
\end{corollary}
\Pf In all cases, $\cH_\lambda$ has several connected components. Thus 
Proposition~\ref{trapping} gives the result. \qed

However one can find non-trapping situations as in Examples~\ref{repulsive-non-trapping} 
and~\ref{1attractive-non-trapping} below. 
\begin{example}\label{repulsive-non-trapping}
Let $N'\in\N^\ast$. Let $V$ be defined on $\R^d\setminus \cS$ by $V(x)=
\sum _{j=1}^{N'}f_j/|x-s_j|$ with
$f_j>0$, for any $1\leq j\leq N'$. It satisfies (\ref{sing-sj}) with $N=0$. 
For $0<\lambda <(\sum _{j=1}^{N'}|s_j|/f_j)^{-1}$, $p$ is non-trapping
at energy $\lambda$. 
\end{example}
\Pf Recall that $a_0(x,\xi)=x\cdot \xi$. For $(x,\xi)\in p^{-1}(\lambda)$, 
\begin{eqnarray*}
\{p,a_0\}(x,\xi)&=&2|\xi |^2\, +\, \sum _{j=1}^{N'}f_jx\cdot
\frac{x-s_j}{|x-s_j|^3}\ =\ |\xi |^2\, +\, \lambda \, +\, \sum _{j=1}^{N'}f_js_j\cdot
\frac{x-s_j}{|x-s_j|^3}\\
&\geq &\lambda \, -\, \sum _{j=1}^{N'}\frac{|s_j|\lambda ^2}{f_j}\ = \ \lambda 
\bigl( 1\, -\, \lambda \sum _{j=1}^{N'}\frac{|s_j|}{f_j}\bigr)\, >\, 0\period
\end{eqnarray*}
Here we used that $0<f_j/|x-s_j|\leq \lambda$, for $(x,\xi)\in
p^{-1}(\lambda)$. Now standard arguments yields the result (see the
proof of Lemma~\ref{escape-infinity}, for instance). \qed 

\begin{example}\label{1attractive-non-trapping}
Let $\lambda, c,\rho >0$ and $W\in C^\infty(\R ^d;\R)$ such that 
\[\forall \alpha\in\N ^d\comma \, \exists C_\alpha >0\, ; \, \forall x\in\R ^d\comma \, 
| \partial _x^\alpha W (x)|\ \leq \ C_\alpha \langle x\rangle ^{-\rho -|\alpha |}\period \]
Let $V\in C^\infty(\R ^d\setminus \{0\};\R)$ defined by $V(x)=-c/|x|+W(x)$. 
Depending on $c$ and $\lambda$, one can find small enough $(C_\alpha)_{|\alpha |\leq 1}$'s 
such that $p$ is non-trapping at energy $\lambda$. 
\end{example}
\Pf The function $\tilde{p}$ defined just before (\ref{cauchy-pb-extended}) takes the 
following form: $\tilde{p}(\pht ;\phz)=\tilde{p}(t,\tau ; z,\zeta )=|\zeta |^2-c+|z|^2
(W\circ \cK (z)-\tau )$. Let $b_0:T^\ast\R\times T^\ast\R ^d\dans \R$ be defined by 
$b_0(\pht ;\phz)=\zeta \cdot z$. Then 
\begin{equation}\label{poisson-upstairs}
\{\tilde{p}\, ,\, b_0\}(\pht ;\phz)\ =\ 2\tilde{p}(\pht ;\phz)\, +\, 2c\, +\, |z|^2\bigl(
4\tau - 4W\circ \cK (z)-z\cdot \nabla _z(W\circ \cK)(z)\bigr)\period 
\end{equation}
Thanks to (\ref{def-KS}), we can choose the $(C_\alpha)_{|\alpha |\leq 1}$'s small enough 
such that, for $\tau=\lambda>0$, the last term in (\ref{poisson-upstairs}) is everywhere 
non-negative. Thus, on $\tilde{p}^{-1}(]-c/2;c/2[)$, $\{\tilde{p}\, ,\, b_0\}\geq c$. 
This implies that, for any solution $s\donne (t(s),\lambda ; z(s),\zeta (s))$ of 
(\ref{cauchy-pb-extended}) leaving in $\tilde{p}^{-1}(0)$, the function $s\donne |z(s)|^2$ 
is strictly convex. It must go to infinity in both time $s$ directions. By (\ref{phi-phitilde}), 
this implies that any broken trajectory $(\phi(t;\ph))_{t\in\R\setminus \J(\ph)}$ 
with $p(\ph)=\lambda$ goes to infinity in both time $t$ directions. \qed 
\begin{remark}\label{rem:5.5}
By inspection of (\ref{poisson-upstairs}) we see that, for a potential
of the form $V(x)=\frac{f(x)}{|x|}+W(x)$ with $f(0)<0$ and meeting (\ref{decroit}), no 
trapping occurs for high enough energies.
\end{remark}

\section{Scattering by a molecular potential.}
\label{sectonex}
\setcounter{equation}{0}

We now show that our analysis can be applied to Example~\ref{ex:mol}.

The potential $x\donne \sum _{j}e_0z_j|x-s_j|^{-1}$ of $P_1(h_0)$ is smooth on
$\hat{M}={\R}^3\backslash\{s_1,\ldots,s_{N'}\}$ and satisfies (\ref{decroit}). 
By local elliptic regularity (see \cite{rs2}, Thm. IX.26), 
the electronic eigenfunctions $\psi _k\in \rL^2(\R^3)$ of $P_1(h_0)$ are smooth on 
$\hat{M}$. Furthermore, they are continuous on ${\R}^3$ (see \cite{cfks}, Thm. 2.4) and 
the corresponding eigenvalues $E_k$ are negative by \cite{fh} (see also \cite{cfks}, Thm. 4.19). 
Using \cite{ag} outside the ball $B:=\{x\in {\R}^3;|x|\leq R_0\}$ (cf. (\ref{decroit})), one can show that 
the $\psi _k$'s decay exponentially. This means, for any $k$, that there exists $c_k,C_k>0$ such that 
\begin{equation}\label{exp-decay}
x\not\in B\ \impl \ |\psi _k(x)|\ \leq \ C_ke^{-c_k|x|}\period 
\end{equation}
By (\ref{decroit}), the result in \cite{ag} 
can be applied to the derivatives of the $\psi _k$ outside $B$. Thus (\ref{exp-decay}) 
holds true for these derivatives with possibly different constants $c_k,C_k$. 
For any $j\in \{1,\ldots ,N'\}$, it turns out that $\psi _{kj}:\R^4\ni z\donne \psi _k(s_j+\cK (z))$, 
with $\cK (z)$ defined in (\ref{def-KS}), is smooth near $z=0$. Indeed,
we can show as in \cite{gk} (see also the proof of
(\ref{tilde-Pj-small}) in Lemma~\ref{tilde-g_n}) that the equation 
$P_1(h_0)\psi _k=E_k\psi _k$ can be lifted to a Schr\"odinger equation in $\R ^4$
with smooth potential solved by the function $\psi _{kj}$. Again, the 
elliptic regularity gives the desired result. Therefore the charge densities $\rho _k:=|\psi _k|^2$ 
are smooth on $\hat{M}$ and continuous on ${\R}^d$. The $\rho _k$ and their derivatives satisfy 
(\ref{exp-decay}). For any $j\in \{1,\ldots ,N'\}$, $\rho
_{kj}:\R^4\ni z\donne \rho _k(s_j+\cK (z))$ is smooth near $z=0$. 
This allows us to obtain the following properties for the $W_k$. 
\begin{proposition}\label{prop-W_k}
Let $k\in\{1,\ldots,K\}$. The potential $W_k$ is smooth on $\hat{M}$, continuous on $\R ^3$, 
and satisfies (\ref{decroit}). For any $j\in \{1,\ldots ,N'\}$, the
function $W _{kj}:\R^4\ni z\donne W _k(s_j+\cK (z))$, 
with $\cK (z)$ defined in (\ref{def-KS}), is smooth near $z=0$.
\end{proposition}
\Pf Since $|\cdot|^{-1}\in \rL ^1(\R ^3)+\rL^\infty (\R ^3)$, $\rho _k\in \rL ^1(\R ^3)$,
and $\rho _k$ is continuous, $W_k$ is well defined and continuous on $\R ^3$. 
Let $j\in \{1,\ldots ,N'\}$, $y\in\hat{M}$, and consider a partition of unity in $\R ^3$ 
of the form $\sum _{j=0}^{N'}\chi _j=1$ with $\chi _j\in C_0^\infty$ and $\chi _j=1$ near $s_j$, 
for $j\geq 1$, and $\chi _0=1$ near $y$. Denoting by $\ast$ the
convolution product, we can write near $y$, for any $k$ and any multiindex
$\alpha\in {\N}^d$,
\begin{equation}\label{deriv-convo}
D_x^\alpha W_k =\ D_x^\alpha (\rho _k \ast
|\cdot|^{-1})\ =\ \sum _{j=1}^{N'}(\rho _k\chi _j)\ast D_x^\alpha
|\cdot|^{-1}\, \, +(D_x^\alpha(\rho _k\chi _0))\ast|\cdot|^{-1} 
\end{equation}
(as distributions). This defines a continuous function near $y$. 
Using the exponential decay of the functions $\rho _k$, we can 
show that $W_k$ satisfies (\ref{decroit}). \\
Let $j\in \{1,\ldots ,N'\}$. We want to show that the function 
$\R ^4\ni z\donne (\rho _k\ast |\cdot|^{-1})({s_j+\cal K}(z))$ 
is a constant times the function $\R ^4\ni z\donne (\rho _{kj}\ast |\cdot|^{-2})(z)$. 
Notice that, for an $f\in C(\R ^4)\cap \rL ^1(\R ^4)$, $f*|\cdot|^{-2}$ 
is a well defined continuous function since $|\cdot|^{-2}\in \rL ^1(\R ^4)+\rL ^\infty (\R ^4)$. 
Now, it is convenient to view $\R ^4$ as the quaternion space $\bH$ and to use the representation 
of $\cK$ on this space (see the appendix). In particular, one can use 
formula (3) from \cite{gk}, saying that for $x:={\cal K}(Y)$, $Y\in \bH$, 
$|Y|^2\,dY=c\cdot dx\,d\theta $
for some constant $c>0$ ($d\theta$ is uniquely defined by
(\ref{local-diffeo}), compare also with the group action (\ref{ga})).
Then, using Lemma  \ref{lemma:fibre} below, we get 
\begin{eqnarray*}
\Big(\rho _{kj}*|\cdot|^{-2}\Big)(Z)&=&\int_{{\R}^4}\frac{\rho _k(s_j+{\cal K}(Y))}{|Y-Z|^2}|Y|^2dY\ =\ c\cdot\int_{{\R}^3\times S^1}\frac{\rho _k(s_j+x)}{|Y(x,\theta)-Z|^2}dx\,d\theta\\
&=&c'\cdot\int_{{\R}^3}\frac{\rho _k(s_j+x)}{|x-{\cal K}(Z)|}dx
\ =\ c'\cdot\int_{{\R}^3}\frac{\rho _k(s_j+x)}{|s_j+x-s_j-{\cal K}(Z)|}dx\\
&=&c'(\rho _k\ast |\cdot|^{-1})(s_j+{\cal K}(Z))\comma 
\end{eqnarray*}
for $Z\in \bH$, and $c'>0$. Now, since $\rho _{kj}$ is smooth near $0$ 
and $|\cdot|^{-1}$ is smooth away from $0$, we can use a formula similar to 
(\ref{deriv-convo}) to show that $\rho _{kj}*|\cdot|^{-2}$ is smooth near $0$. \qed

\begin{lemma}\label{lemma:fibre}
For $X,Z\in\bH$ with ${\cal K}(Z)\neq {\cal K}(X)$
\[\int_0^{2\pi} |\exp(I_1\theta)Z-X|^{-2}\,d\theta
= 2\pi \cdot |{\cal K}(Z)-{\cal K}(X)|^{-1}\period \]
\end{lemma}
\Pf 
Assuming the condition, both sides are well-defined. Then, using the
definition of the real part of a quaternion (see appendix)
\beqno
\lefteqn{\int_0^{2\pi} |\exp(I_1\theta)Z-X|^{-2}\,d\theta}\\
&=& \int_0^{2\pi} \l(|Z|^2+|X|^2-2{\rm
Re}((\cos(\theta)+I_1\sin(\theta))ZX^\ast)\ri)^{-1}\,d\theta\\
&=& \int_0^{2\pi} \l(|Z|^2+|X|^2-2({\rm
Re}(ZX^\ast)\cos(\theta)+{\rm
Re}(I_1ZX^\ast)\sin(\theta))\ri)^{-1}\,d\theta\\
&=& \int_0^{2\pi} \l(|Z|^2+|X|^2-2\sqrt{
({\rm Re}(ZX^\ast))^2+({\rm Re}(I_1ZX^\ast))^2}\cos(\psi)\ri)^{-1}\,d\psi\\
&=& 2\pi \cdot
\big((|Z|^2+|X|^2)^2-
4({\rm Re}(ZX^\ast))^2-4({\rm Re}(I_1ZX^\ast))^2\big)^{-1/2}\\
&=& 2\pi \cdot |Z^\ast I_1Z-X^\ast I_1X|^{-1}\ = \ 
2\pi \cdot |{\cal K}(Z)-{\cal K}(X)|^{-1}\comma 
\eeqno

the last two equations being due to (\ref{identity}) and (\ref{Hopf:map}). \qed

Now we are able to explain why the proof of our results can be adapted to 
treat the potential $V$ defined in (\ref{VW}). In the proof of the necessity 
of the non-trapping condition in Section~\ref{resolv-esti-non-trapp}, 
the results away from the singularities work since $V$ satisfies (\ref{decroit}). 
Since the $W _{kj}$ are smooth near $0$, the results in \cite{gk} (see 
Lemmata~\ref{sing-evolution-u_h} and~\ref{time-space-localization}) are 
still valid. Since each $W_k$ is bounded, it is small compared to a repulsive 
potential $+|\cdot -s_j|^{-1}$ near the corresponding repulsive singularity $s_j$. 
So Section~\ref{repulsive-case} is also valid. In the proof of the converse 
in Section~\ref{semicl-trapp}, the results away from the singularities hold true 
since (\ref{decroit}) is still valid. The fact that the $W_k$ is small compared to the size of 
a singular potential $\pm |\cdot -s_j|^{-1}$ near the corresponding singularity $s_j$ 
explains why Section~\ref{repuls-sing} works and also the validity of (\ref{leading-term}). 
The fact that the $W _{kj}$ are smooth near $0$, ensures that 
Lemma~\ref{tilde-g_n} still works.

\appendix
\section{The Hopf map.}
We use the following notation for the quaternion algebra over ${\R}$:
\[\bH := \l\{\l.\bsm w_1 & -w_2\\ \bar{w}_2 &
\bar{w}_1\esm 
\ri| w_1,w_2\in {\mathbb C}\ri\} \cong{\R}^4\]  
with matrix multiplication, and basis 
\beqno
(I_0,I_1,I_2,I_3)&:=&\l(
\bsm 1& 0 \\ 0& 1 \esm , 
\bsm i& 0 \\ 0& -i\esm ,
\bsm 0& 1 \\-1& 0 \esm ,
\bsm 0& -i \\ -i& 0 \esm 
\ri).
\eeqno 
The direct sum decomposition 
$\bH= {\R}\cdot \idty \oplus {\rm Im} \bH$
with 
\beqno
{{\rm Im}\bH} &:=& \{Z\in\bH\mid Z^2 = \lambda \cdot \idty\mbox{ with }
\lambda \leq 0 \}={\rm Span}_{\R} (I_1,I_2,I_3)
\eeqno
into real and imaginary space
is orthogonal w.r.t.\ the inner product
\[\bH\times\bH \ar{\R} \quad,\qquad \LA X,Y\RA := \eh \tr(XY^\ast),\]
$X\mapsto X^\ast:=\bar{X}^t$ being the conjugation.
The norm $|X|:=\LA X,X\RA^{\eh}$ is multiplicative:
\[|XY| = |X|\,|Y|\qquad (X,Y\in\bH).\]
The {\em real part} of a quaternion equals ${\rm Re}(X):= \eh
\tr(X)$.\\
See, e.g., \cite{EHKKMNPRE} for more information on $\bH$. The {\em Hopf map} equals 
\beq
{\cal K}:\bH\ar{{\rm Im}\bH}\quad,\quad {\cal K}(Z):=Z^\ast I_1 Z=
i\bem{cc} 
w_1\bar{w}_1-w_2\bar{w}_2 & -2\bar{w}_1 w_2 \\ 
-2w_1 \bar{w}_2      & w_2\bar{w}_2-w_1\bar{w}_1\eem 
\Leq{Hopf:map}
which is a surjection ${\R}^4\ar{\R}^3$ whose preimages are the orbits of the
isometric group action
\beq
\alpha_0:S^1\ar{\rm Aut}(\bH),\qquad\alpha_0(\theta)(Z):=
\exp(\theta I_1)Z.  
\Leq{ga}
This action is free on $\bH\setminus \{0\}$.
We call ${\cal K}$ the Hopf map,  since its restriction to $S^3$ is the Hopf fibration 
$S^3\dans S^2$ with fibre $S^1$. 

Writing $w_1:=z_0+i z_3$, $w_2:=z_2+i z_1$
we get formula (\ref{def-KS}) in the basis 
$(I_1,I_2,I_3)$ of ${\rm Im}\bH$. Finally we prove the formula 
\begin{equation}\label{identity}
|Z^\ast I_1 Z-X^\ast I_1 X|\ =\ \sqrt{(|Z|^2+|X|^2)^2-
4\bigl(({\rm Re}(ZX^\ast))^2+({\rm Re}(I_1ZX^\ast))^2\bigr)}
\end{equation}
used in Section~\ref{sectonex}. \\
Notice that, for all $A,B\in\bH$, ${\rm Re}(I_kA^\ast)=
-{\rm Re}(I_kA)$, ${\rm Re}(A^\ast)={\rm Re}(A)$, ${\rm Re}(A^\ast A)=|A|^2$, and 
\begin{equation}\label{rere}
{\rm Re}(AB) =
{\rm Re}(A){\rm Re}(B)- \sum_{k=1}^3{\rm Re}(I_kA){\rm Re}(I_kB)\period 
\end{equation}
Setting $A:=I_1ZX^\ast$ and $B:= I_1XZ^\ast$ in (\ref{rere}), we get
\[{\rm Re}\big((I_1 ZX^\ast)(I_1 XZ^\ast)\big)=
-({\rm Re}(XZ^\ast))^2 - ({\rm Re}(I_1XZ^\ast))^2 + ({\rm Re}(I_2XZ^\ast))^2 + ({\rm Re}(I_3XZ^\ast))^2\]

Similarly it follows from (\ref{rere}) that 
$|A|^2=\sum_{k=0}^3({\rm Re}(I_k A))^2$, so that for 
$A:=ZX^\ast$
\[|Z|^2\,|X|^2=|A|^2= 
({\rm Re}(ZX^\ast))^2  + ({\rm Re}(I_1ZX^\ast))^2 +
({\rm Re}(I_2XZ^\ast))^2+ ({\rm Re}(I_3XZ^\ast))^2.\]
So
\beqno
\lefteqn{|Z^\ast I_1Z-X^\ast I_1X|^2=(Z^\ast I_1Z-X^\ast I_1X)(-Z^\ast I_1Z+X^\ast I_1X)}\\
&=&
Z^\ast I_1(-|Z|^2)I_1Z+X^\ast I_1(-|X|^2)I_1X+(Z^\ast I_1ZX^\ast I_1X)+(Z^\ast I_1ZX^\ast I_1X)^\ast\\
&=& |Z|^4+|X|^4+2{\rm Re}\big((I_1ZX^\ast)(I_1XZ^\ast)\big)\\
&=&
|Z|^4+|X|^4+2(-({\rm Re}(XZ^\ast))^2 - ({\rm Re}(I_1XZ^\ast))^2  + ({\rm Re}(I_2 XZ^\ast))^2 + ({\rm Re}(I_3 XZ^\ast))^2
\\
&=&(|Z|^2+|X|^2)^2-
4\l(({\rm Re}(ZX^\ast))^2+({\rm Re}(I_1ZX^\ast))^2\ri)
\eeqno
This proves the claim.\hfill\qed

{\bf Acknowledgements.}
The authors would like to thank the anonymous referee
for numerous valuable remarks and suggestions, that helped to improve the clarity of this text,
as well as for pointing out reference \cite{w3}.


%
%
\end{document}